\documentclass[11pt]{amsart}
\usepackage{amssymb,mathrsfs}

\usepackage{amssymb,mathrsfs}
\usepackage{tikz-cd}
\usepackage{dynkin-diagrams}

\usepackage{fancyhdr}

\begin{document}

\title{Stable Automorphic Forms for Semisimple Groups }

\author{Jae-Hyun Yang}

\address{Yang Institute for Advanced Study
\newline\indent
130 Mokdongseo-ro Yangcheon-gu
\newline\indent
Seoul 07989, Korea
\vskip 2mm
and
\vskip 2mm
Department of Mathematics
\newline\indent
Inha University
\newline\indent
Incheon 22212, Korea}

\email{jhyang@inha.ac.kr\ \ or\ \ jhyang8357@gmail.com}

\newtheorem{theorem}{Theorem}[section]
\newtheorem{lemma}[theorem]{Lemma}
\newtheorem{corollary}[theorem]{Corollary}
\newtheorem{proposition}[theorem]{Proposition}
\newtheorem{remark}[theorem]{Remark}
\newtheorem{definition}[theorem]{Definition}
\newtheorem{conjecture}[theorem]{Conjecture}
\newtheorem{example}[theorem]{Example}
\newtheorem{exercise}[theorem]{Exercise}
\newtheorem{problem}[theorem]{Problem}
\newtheorem{question}[theorem]{Question}

\renewcommand{\theequation}{\thesection.\arabic{equation}}
\renewcommand{\thetheorem}{\thesection.\arabic{theorem}}
\renewcommand{\thelemma}{\thesection.\arabic{lemma}}

\newcommand{\BR}{\mathbb R}
\newcommand{\BQ}{\mathbb Q}
\newcommand{\BT}{\mathbb T}
\newcommand{\BM}{\mathbb M}
\newcommand{\bn}{\bf n}
\def\charf {\mbox{{\text 1}\kern-.24em {\text l}}}
\newcommand{\BC}{\mathbb C}
\newcommand{\BZ}{\mathbb Z}

\thanks{2020 {\it Mathematics\ Subject\ Classification.} Primary 11F46, 11F55, 11G10, 11G15, 32N10.
\endgraf
{\it Keywords\ and\ phrases.} stable automorphic forms, universal arithmetic varieties, stable Schottky-Siegel \\
\indent forms.}

\begin{abstract}
In this paper, we introduce the concept of stable automorphic forms for semisimple algebraic groups and use the stability of automorphic forms to study the geometry of infinite dimensional arithmetic quotients.
\end{abstract}

\maketitle

\newcommand\tr{\triangleright}
\newcommand\al{\alpha}
\newcommand\be{\beta}
\newcommand\g{\gamma}
\newcommand\gh{\Cal G^J}
\newcommand\G{\Gamma}
\newcommand\de{\delta}
\newcommand\e{\epsilon}
\newcommand\lam{\lambda}
\newcommand\z{\zeta}
\newcommand\vth{\vartheta}
\newcommand\vp{\varphi}
\newcommand\om{\omega}
\newcommand\p{\pi}
\newcommand\la{\lambda}
\newcommand\lb{\lbrace}
\newcommand\lk{\lbrack}
\newcommand\rb{\rbrace}
\newcommand\rk{\rbrack}
\newcommand\s{\sigma}
\newcommand\w{\wedge}
\newcommand\fgj{{\frak g}^J}
\newcommand\lrt{\longrightarrow}
\newcommand\lmt{\longmapsto}
\newcommand\lmk{(\lambda,\mu,\kappa)}
\newcommand\Om{\Omega}
\newcommand\ka{\kappa}
\newcommand\ba{\backslash}
\newcommand\ph{\phi}
\newcommand\M{{\Cal M}}
\newcommand\bA{\bold A}
\newcommand\bH{\bold H}
\newcommand\D{\Delta}

\newcommand\Hom{\text{Hom}}
\newcommand\cP{\Cal P}

\newcommand\cH{\Cal H}

\newcommand\pa{\partial}

\newcommand\pis{\pi i \sigma}
\newcommand\sd{\,\,{\vartriangleright}\kern -1.0ex{<}\,}
\newcommand\wt{\widetilde}
\newcommand\fg{\frak g}
\newcommand\fk{\frak k}
\newcommand\fp{\frak p}
\newcommand\fs{\frak s}
\newcommand\fh{\frak h}
\newcommand\Cal{\mathcal}

\newcommand\fn{{\frak n}}
\newcommand\fa{{\frak a}}
\newcommand\fm{{\frak m}}
\newcommand\fq{{\frak q}}
\newcommand\CP{{\mathscr P}_n}
\newcommand\Hnm{{\mathbb H}_n \times {\mathbb C}^{(m,n)}}
\newcommand\BD{\mathbb D}
\newcommand\BH{\mathbb H}
\newcommand\CCF{{\mathcal F}_n}
\newcommand\CM{{\mathcal M}}
\newcommand\Gnm{\Gamma_{n,m}}
\newcommand\Cmn{{\mathbb C}^{(m,n)}}
\newcommand\Yd{{{\partial}\over {\partial Y}}}
\newcommand\Vd{{{\partial}\over {\partial V}}}

\newcommand\Ys{Y^{\ast}}
\newcommand\Vs{V^{\ast}}
\newcommand\LO{L_{\Omega}}
\newcommand\fac{{\frak a}_{\mathbb C}^{\ast}}

\vskip 1mm

\vskip 7mm

\centerline{\bf Table of Contents}

\vskip 0.75cm $ \quad\qquad\textsf{\large \ 1.
Introduction}$\vskip 0.03251cm

$\quad\qquad \textsf{\large\ 2. Stable functions on infinite dimensional varieties}$
\vskip 0.043251cm

$ \quad\qquad  \textsf{\large\ 3. Stable automorphic forms for a semisimple algebraic group }$
\vskip 0.04321cm

$ \quad\qquad \textsf{\large\ 4. Examples of stable automorphic forms}$
\vskip 0.043251cm
$ \quad\qquad\quad
\textsf{\large\ 4.1. Stable automorphic forms for $Sp(\infty,\BR)$}$
\vskip 0.043251cm
$ \quad\qquad\quad
\textsf{\large\ 4.2. Stable automorphic forms for $SL(\infty,\BR)$}$
\vskip 0.043251cm

$ \quad\qquad  \textsf{\large\ 5. Applications of the stability to geometry}$
\vskip 0.043251cm

$ \quad\qquad\quad
\textsf{\large\ 5.1. The universal moduli space of abelian varieties}$
\vskip 0.043251cm

$ \quad\qquad\quad
\textsf{\large\ 5.2. The universal moduli space of curves}$
\vskip 0.043251cm

$ \quad\qquad\quad
\textsf{\large\ 5.3. The universal moduli space of polarized real tori}$
\vskip 0.0437251cm


$ \quad\qquad\ \textsf{\large References }$

\vskip 10mm


\begin{section}{{\large\bf Introduction}}
\setcounter{equation}{0}
\vskip 0.1cm

Originally the notion of stable automorphic forms was introduced
in the symplectic group by E. Freitag \cite{Fr2} in 1977. Those automorphic forms were
called {\it stable\ modular\ forms} by Freitag. He proved that the set of all stable modular forms is a polynomial ring in a countably infinite set of indeterminates over the field $\BC$ of complex numbers
(cf. \cite[Theorem 2.5,\,p.\,204]{Fr2}). Thereafter R. Weissauer
investigated stable modular forms in the sense of Freitag intensively
for the study of Eisenstein series \cite{We}.
In 2014, Codogni and Shepherd-Barron \cite{Cod-Sh} proved that there do not exist nontrivial {\it stable}
{\it Schottky}-{\it Siegel\ modular\ forms} for the Jacobian locus. In 2016,
Codogni \cite{Cod} proved that there exist nontrivial {\it stable}
{\it Schottky}-{\it Siegel\ modular\ forms} for the hyperelliptic locus.
\vskip 3mm

In this article, we will deal with the case of stable automorphic forms only for semisimple algebraic groups. The motivation of introducing the notion of stable automorphic forms is to investigate
the geometric properties of finite or infinite dimensional arithmetic quotients
associated with those automorphic forms.

\vskip 2mm
The purpose of this article is to generalize the concept of stable modular forms to that
of stable automorphic forms for semisimple algebraic groups
and apply the stability of automorphic forms to the study of the universal moduli
space of abelian varieties, the universal moduli space of curves and the universal moduli
space of polarized real tori.
\vskip 2mm
This paper is organized as follows. In section 2, we review the notion of
infinite dimensional algebraic varieties due to I.\,R. Shafarevich\,\cite{Sh1,Sh2,K}.
We introduce the notion of stable functions. In section 3, we introduce the
notion of stable automorphic forms for semisimple algebraic groups.
In the case of a finite dimensional semisimple algebraic group, we follow the definition of
automorphic forms given by Harish-Chandra\,(cf.\,\cite{HC3}) and Borel\,(cf.\,\cite{Bo, B-J}).
In section 4, as examples, we consider stable automorphic forms for both an infinite dimensional symplectic group $Sp(\infty,\BR)$ and an infinite dimensional special linear group $SL(\infty,\BR)$. In the final section, using the stability of automorphic forms for
$Sp(\infty,\BR)$ and $SL(\infty,\BR)$, we characterize the so-called {\it universal\,(or\ stable)\ Satake\ compactifications} and investigate their geometry. We deal with the universal moduli space of abelian varieties, the universal moduli space of curves and the universal moduli space of polarized real tori.

\vskip 3mm
\noindent
{\bf Notations.} We denote by $\BZ,\ \BR$ and $\BC$ the ring of integers, the field of real numbers and the field of complex numbers respectively.
$\BZ^+$ and $\BZ_+$ denote the set of all positive integers and the set of
all nonnegative integers respectively.
The symbol ``:='' means that the expression on the right is the definition
of that on the left. $F^{(k,l)}$ denotes the set of all $k\times l$ matrices
with entries in a commutative ring $F$. For a square matrix $A,\ {\rm Tr}(A)$
denotes the trace of $A$. For any $M\in F^{(k,l)},\ ^tM$ denotes the
transpose of $M$. For $A\in F^{(k,l)}$ and $B\in F^{(k,k)},$ we set
$B[A]=\,^tABA$\,(Siegel's notation). $I_n$ denotes the identity matrix of degree $n$.
For a complex matrix $A$, ${\overline A}$ denotes the complex {\it conjugate} of $A$.
${\rm diag}(a_1,a_2,\cdots,a_n)$ denotes the $n\times n$ diagonal matrix with diagonal entries
$a_1,\cdots,a_n$. For a smooth manifold, we denote by $C_c (X)$ (resp. $C_c^{\infty}(X)$ the algebra of all continuous (resp. infinitely differentiable) functions on $X$ with compact support.
\begin{equation*}
J_n=\begin{pmatrix} 0&I_n\\
                   -I_n&0\end{pmatrix}
\end{equation*}
denotes the symplectic matrix of degree $2n$.
\begin{equation*}
{\mathbb H}_n=\,\{\,\Omega\in \BC^{(n,n)}\,|\ \Omega=\,^t\Omega,\ \ \ \text{Im}\,\Omega>0\,\}
\end{equation*}
denotes the Siegel upper half plane of degree $n$.
\begin{equation*}
  Sp(2n,\BR)=\{ M\in \BR^{(2n,2n)}\,|\ {}^tM J_nM=J_n\,\}
\end{equation*}
denotes the symplectic group of degree $n$ and
\begin{equation*}
  \G_n=\{ \gamma\in \BZ^{(2n,2n)}\,|\ {}^t\gamma J_n\gamma=J_n\,\}\subset Sp(2n,\BR)
\end{equation*}
denotes the Siegel modular group of degree $n$.
For a positive integer $n$, we denote
$$\mathscr{P}_n=\{\,Y\in\BR^{(n,n)}\,|\ Y=\,{}^tY>0\,\}$$
and
$$ \mathbb{X}_n=\{\,Y\in\BR^{(n,n)}\,|\ Y=\,{}^tY>0,\ \det\,(Y)=1\,\}.$$
We denote $\mathfrak{X}_n=SL(n,\BZ)\ba \mathbb{X}_n.$

\end{section}

\vskip 7mm

\begin{section}{{\large\bf Stable functions on infinite dimensional varieties}}
\setcounter{equation}{0}
\vskip 0.2cm
First we review the notion of infinite dimensional algebraic groups
due to I.\,R. Shafarevich\,(\,cf.\,\cite{Sh1, Sh2, K}).

\noindent
\begin{definition}\label{def:2.1}
By an {\sf infinite dimensional algebraic variety} over a field $k$
we mean the inductive limit $X$ of a directed system $(\,X_i,\,f_{ij}\,)$
of finite dimensional algebraic varieties over
the field $k$, where $f_{ij}:X_i\lrt X_j\,(\,i<j\,)$ are closed
embeddings. We write
\begin{equation*}
 X:=\,
 \lim_{\begin{subarray}{c} \longrightarrow\\ ^i \end{subarray}}X_i.
\end{equation*}
\end{definition}

\vskip 2mm
Throughout this paper, we shall consider only the case where the set of indices is the set $\BZ^+$ of all positive integers. Each of the
$X_i$ will be considered to be equipped with its Zariski topology and
we endow $X$ with the topology of the inductive limit where a set
$Z\subset X$ is closed if and only if its preimage in each $X_i$
is closed. In particular, each $X_i$ is closed in $X.$

\vskip 1mm
\begin{definition}\label{def:2.2}
A continuous mapping $f:X\lrt Y$ of two infinite
dimensional algebraic varieties is called a {\sf morphism} if for
any $X_i$ in the system $(X_i)$ defining $X$, there exist at least one $Y_j$ in the system $(Y_j)$ defining $Y$ such that
$f(X_i)\subset Y_j$ and the restriction $f:X_i\lrt Y_j$ is a
morphism of finite dimensional algebraic varieties. Irreducibility and connectedness of an infinite dimensional algebraic variety are
defined as irreducibilty and connectedness of the corresponding
topological space.
\end{definition}

\vskip 1mm
\begin{definition}\label{def:2.3}
An infinite dimensional algebraic variety $G$
with a group structure is called an {\sf infinite\ dimensional\
algebraic\ group} if the inverse mapping $x\longmapsto x^{-1}$ and
the multiplication $(x,y)\longmapsto xy$ are morphisms for all $x,y\in G$.
\end{definition}

\vskip 2mm
In a similar way, we may define the notions of infinite dimensional
smooth manifolds, infinite dimensional complex manifolds,
infinite dimensional real or complex Lie groups and so on with a
usual topology and suitable morphisms.

\vskip 2mm
Let $X$ be an infinite dimensional space with its directed system
$(\,X_i,\,f_{ij}\,)$. Let $V$ be a fixed finite dimensional complex
vector space. We assume that

\vskip 2mm
(I) to each $X_i$ there is given the vector space $C_i$ of functions
on $X_i$ with values in $V$
\par \ \ \ \,
and that
\vskip 1mm
(II) there is given an inverse system $(\,C_i,\,\Phi_{ij}\,)$ of
linear maps $\Phi_{ij}:C_j\lrt C_i\,(\,i<j\,)$ \\
\indent \ \ \ \ \ \,such that
$$\Phi_{ik}\,=\,\Phi_{ij}\circ \Phi_{jk}\ \ \ \ \text{for\ all}\
i<j<k.$$
Now we let
\begin{eqnarray*}
C:&=&\,\lim_{\begin{subarray}{c} \longleftarrow\\ ^i \end{subarray}}
(\,C_i,\,\Phi_{ij}\,)\\
&=&\,\left\{\,(f_k)\in \prod_{i\in \BZ^+}\,C_i\,\vert\ \,
\Phi_{ij}(f_j)=f_i\
\text{for\ all}\ i<j\ \right\}
\end{eqnarray*}
be the inverse limit of the system $(\,C_i,\,\Phi_{ij}\,).$
Elements of $C$ are called {\sf stable\ functions}.
\vskip 2mm
\begin{example}\label{exa:2.4}
Let
\begin{equation*}
  \BR^{\infty}:=\,
 \lim_{\begin{subarray}{c} \longrightarrow\\ ^i \end{subarray}}\BR^i
\end{equation*}
be the infinite dimensional Euclidean space, where $\BR^i$ is the real
Euclidean space of dimension $i$. For a positive integer $i$, we let
$C_c(\BR^i)$ be the vector space of all real-valued continuous functions
on $\BR^i$ with compact support. For any two positive integers $i,j$ with
$i<j$, we define
$$\psi_{ij}:C_c(\BR^j)\lrt C_c(\BR^i)$$
by
\begin{equation*}
  \psi_{ij}(f)(x):=\,f(x,0_{j-i}),\quad x\in\BR^i,\ 0_{j-i}=(0,\cdots,0)\in\BR^{j-i}.
\end{equation*}
We take the inverse limit
\begin{equation*}
\mathscr{C}=\,\lim_{\begin{subarray}{c} \longleftarrow\\ ^i \end{subarray}}
(\,C_c(\BR^i),\,\psi_{ij}\,)
\end{equation*}
of the system $(\,C_c(\BR^i),\,\psi_{ij}\,)$. Elements of $\mathscr{C}$ are stable
continuous functions with compact support.
\end{example}

\vskip 1mm
\begin{example}\label{exa:2.5}
For any two nonnegative integers $k,l\in \BZ_+$ with
$k<l,$ we define the mapping $\varphi_{kl}:\BH_k\lrt \BH_l$ by
\begin{equation*}
\varphi_{kl}(Z):=\,\begin{pmatrix} Z & 0 \\ 0 & iI_{l-k}
\end{pmatrix},\quad Z\in \BH_k.
\end{equation*}
We recall that $\BH_n$ denotes the Siegel upper half plane of degree $n$ (see Notations).
Then the image $\varphi_{kl}(\BH_k)$ is a totally geodesic subspace of $\BH_l$. We let
\begin{equation*}
\BH_\infty=\,\varinjlim_k \BH_k
\end{equation*}
be the inductive limit of the direct system $(\BH_k,\,\varphi_{kl}).$ $\BH_\infty$ can be described explicitly as follows:
\begin{equation*}
  \left\{ \begin{pmatrix}
            Z & 0 \\
            0 & iI_\infty
          \end{pmatrix} \bigm|\ \, Z\in \BH_k
          \ {\rm for\ some}\ k\geq 1 \right\}.
\end{equation*}
We can show that $\BH_\infty$ is an infinite
dimensional smooth Hermitian symmetric manifold locally closed on $\BC^\infty$, the complex vector space of finite sequences with the finite topology\,(cf.\,\cite{Gl, H}). $\BH_\infty$ has an invariant Riemannian metric which induces the normalized Riemannian metric on each embedded interior subspace $\BH_k$ in $\BH_\infty$.
\vskip 2mm
The symplectic group $Sp(2n,\BR)$ acts on $\BH_n$ transitively by
\begin{equation}\label{(2.1)}
  g\cdot\Omega:=(A\Omega +B)(C\Omega +D)^{-1},
\end{equation}
where $g=
\begin{pmatrix}
  A & B \\
  C & D
\end{pmatrix}\in Sp(2n,\BR)$ and $\Omega\in\BH_n.$ For a fixed nonnegative integer $k$,
a holomorphic function $f:\BH_n\lrt \BC$ is called a {\sf Siegel\ modular\ form} of {\it weight} $k$
if it satisfies the following condtions\,:
\vskip 2mm
\noindent
\ \ {\rm (SM1)} $f(\g\cdot \Omega)=\det (C\Omega+D)^k f(\Omega)$\ \ for all
$\g=\begin{pmatrix}
  A & B \\
  C & D
\end{pmatrix}\in \G_n$ and $\Omega\in\BH_n.$
\vskip 2mm
\noindent
\ \ {\rm (SM2)} If $n=1$, $f$ requires a cuspidal condition, that is, $f$ is
bounded in any domain \\
\indent \ \ \ \ \ \ \ \
$y\geq y_0>0$ and vanishes at the cusp.

\vskip 2mm
\noindent
Here $\G_n$ is the Siegel modular group of degree $n$ (See Notations).
We denote $[\G_n,k]$ the vector space of all Siegel modular forms of weight $k$.
For any two positive integers $i,j$ with $i<j$, we recall the well-known Siegel operator
\begin{equation*}
  \varphi_{ij}:[\G_j,k]\lrt [\G_i,k]
\end{equation*}
defined by
\begin{equation}\label{(2.2)}
  \varphi_{ij}(f)(\Omega):=\,
 \lim_{t\lrt \infty} f
 \begin{pmatrix}
  \Omega & 0 \\
  0 & \sqrt{-1}\,t\,I_{j-i}
\end{pmatrix},\quad \Omega\in\BH_i.
\end{equation}
It is well known that $\varphi_{ij}$ is a well-defined linear map (cf.\,\cite{Fr3}).
For a fixed nonnegative integer $k$, we take the inverse limit
\begin{equation*}
[\G_\infty,k]:= \,\lim_{\begin{subarray}{c} \longleftarrow\\ ^i \end{subarray}}
(\,[\G_i,k],\,\varphi_{ij}\,)
\end{equation*}
of the system $(\,[\G_i,k],\,\varphi_{ij}\,)$. Elements of $[\G_\infty,k]$ are called
{\sf stable\ modular\ forms}.
\end{example}

\end{section}

\vskip 7mm


\begin{section}
{{\large\bf Stable automorphic forms for a semisimple algebraic group}}
\setcounter{equation}{0}

\vskip 2mm
Before we introduce the notion of stable automorphic forms for a semisimple algebraic group,
we recall automorphic forms on a semisimple algebraic group of
finite dimesion\,(cf.\,\cite{Bo,B-J,HC3}).
\vskip 2mm
Let $G$ be a finite dimensional semisimple algebraic group with a maximal compact subgroup $K$.
Let $\rho:K\lrt GL(V)$ be a given representation of $K$ on a finite dimensional complex vector
space $V$. Let $\G$ be an arithmetic subgroup of $G$. A smooth vector-valued function
$f:G\lrt V$ is called an {\sf automorphic form} of type $(\rho,\G)$ if it satisfies the
following conditions (AF1)--(AF3)\,:
\vskip 1mm
(AF1) $f(\g gk)=\rho(k)^{-1}f(g)$ \ for all $\g\in \G,\ g\in G$ and
$k\in K.$
\vskip 2mm
(AF2) $f$ is $Z(\fg)$-finite.
\vskip 2mm
(AF3) $f$ satisfies a suitable growth condition (cf.\,\cite[pp.\,190--191]{B-J}).

\vskip 2mm\noindent
Here $Z(\fg)$ denotes the center of the universal enveloping algebra $U(\fg)$ of
the Lie algebra $\fg$ of $G.$
\vskip 3mm
\begin{theorem}\label{thm:3.1}
 The vector space of all automorphic forms of type $(\rho,\G)$ is finite dimensional.
\end{theorem}
\noindent
{\it Proof.} The proof was done by Harish-Chandra\,\cite{HC3}. \hfill $\square$

\vskip 2mm
Let $G_\infty$ be an infinite dimensional semisimple algebraic group with its inductive system $(\,G_i,\,\phi_{ij}\,)$ of finite dimensional semisimple algebraic groups $G_i$ and the group monomorphisms $\phi_{ij}:G_i\lrt G_j\,(\,i<j\,).$
{\it We\ fix\ a\ finite\ dimensional\ complex\ vector\ space} $V$.

\vskip 3mm
We now assume that
\vskip 2mm
\ \ {\bf (I)} there is given a sequence $(K_i)$ of comapct subgroups such that each $K_i$ is a maximal\\
\indent\ \ \ \ \ \ \ compact subgroup of $G_i$ and $\phi_{ij}(K_i)
\subset K_j$ for all $i<j.$

\vskip 2mm
\ {\bf (II)} there is given a sequence $(\G_i)$ such that each $\G_i$ is an arithmetic subgroup of $G_i$ and
\\
\indent\ \ \ \ \ \ \ $\phi_{ij}(\G_i)\subset \G_j$ for all
$i<j.$

\vskip 2mm
{\bf (III)} there is given a sequence $(\rho_i)$ such that each $\rho_i$ is a representation of $K_i$ on $V$ \\
\indent\ \ \ \ \ \ \ \,
compatible with the morphisms
$\phi_{ij}$, that is, if $i<j,$ then $\rho_j(\phi_{ij}(k))=\rho_i(k)$
\\
\indent\ \ \ \ \ \ \ \,
for all $k\in K_i.$

\vskip 2mm
For each positive integer $i\in \BZ^+,$ we let ${\bf A}(\G_i,\rho_i)$ be the complex vector space
of all automorphic forms of type $(\rho_i,\G_i).$ According to the definition,
we see that if $f\in {\bf A}(\G_i,\rho_i),$ then $f$ satisfies the following conditions
${\rm(AF1)_i}$--${\rm(AF3)_i}$\,:
\vskip 2mm
${\rm(AF1)_i} \ \,f(\g gk)=\rho_i(k)^{-1}f(g)$ \ for all $k\in K_i,\ g\in G_i$ and
$\g\in \G_i.$
\vskip 2mm
${\rm(AF2)_i} \ \,f$ is $Z(\fg_i)$-finite.
\vskip 2mm
${\rm(AF3)_i} \ \,f$ satisfies a suitable growth condition.

\vskip 2mm\noindent
Here $Z(\fg_i)$ denotes the center of the universal enveloping algebra $U(\fg_i)$ of the Lie algebra $\fg_i$ of $G_i.$

\vskip 2mm
We also assume that

\vskip 2mm

{\bf (IV)} there is a sequence $\{ A_i\}$ and $\{ L_{ij}\}$ of linear maps
$$L_{ij}\,:\,A_j\lrt A_i,\ \ \ i<j$$
\\
\indent\ \ \ \ \ \ \ \,
satisfying the conditions
$$L_{ik}\,=\,L_{ij}\circ L_{jk}\ \ \ \text{for\ all}\ i<j<k,$$
where $A_i$ is a subspace of ${\bf A}(\G_i,\rho_i), \ i=1,2,\cdots.$
Elements of the inverse limit of the inverse system $(A_i,L_{ij})$
\begin{equation*}
A_\infty:=\,\lim_{\begin{subarray}{c} \longleftarrow\\ ^i \end{subarray}}
A_i
\end{equation*}
are called {\sf stable\ automorphic\ forms} for an infinite dimensional semisimple algebraic group $G_\infty.$ If there is no confusion, we briefly say stable automorphic forms.

\vskip 2mm
We put
\begin{equation*}
K_\infty:=\,\varinjlim_i K_i,\qquad
\G_\infty:=\,\varinjlim_i \G_i
\ \ \ \text{and}\ \ \ \rho_\infty:=\,\varinjlim_i \rho_i.
\end{equation*}
We call $\rho_\infty$ a {\sf stable\ representation}
of $K_\infty$ or simply a {\sf stable\ representation}. It is easy to see that a stable automorphic form $f$ in $A_\infty$ satisfies the following conditions (S1)--(S3)\,:

\vskip 2mm
\noindent
(SAF1) \ $f(\g gk)\,=\,\rho_\infty (k)^{-1}
f(g)$ \ for all $k\in K_\infty,\ g\in G_\infty$ and $\g\in \G_\infty.$
\vskip 2mm
\noindent
(SAF2) \ $f$ is $Z(\fg_\infty)$-finite.
\vskip 2mm
\noindent
(SAF3) \ $f$ satisfies a suitable growth condition.
\vskip 2mm
\noindent
Here
$Z(\fg_\infty)$ denotes the center of the universal enveloping algebra $U(\fg_\infty)$ of
the Lie algebra $\fg_\infty$ of an infinite dimensional semisimple algebraic group $G_\infty.$

\begin{problem}\label{prob:3.2}
Investigate the structure of $A_\infty$.
\end{problem}

\vskip 10mm

\end{section}


\begin{section}
{{\large\bf Examples of stable automorphic forms}}
\setcounter{equation}{0}

\vskip 2mm
In this section, we give two examples of stable automorphic forms for both $Sp(\infty,\BR)$ and $SL(\infty,\BR)$.

\newcommand\PZB{ {{\partial}\over {\partial{\overline Z}}} }
\newcommand\PZ{ {{\partial}\over{\partial Z}} }
\vskip 3mm
\noindent
{\bf Example A.\ Stable automorphic forms for $Sp(\infty,\BR)$}
\vskip 3mm
First of all, we provide some geometric properties on the Einstein-K{\"a}hler Hermitian
symmetric manifold $Sp(2n,\BR)/U(n)\cong \BH_n$ that is important geometrically and
number theoretically.
\vskip 2mm
We let $G:=Sp(2n,\BR)$ and $K=U(n).$ We recall that $G$ acts on ${\mathbb H}_n$
transitively via the formula (2.3).
The stabilizer of the action (2.3) at $iI_n$
is
\begin{equation*}
  \left\{ \begin{pmatrix} \,A & B \\ -B & A \end{pmatrix} \Big| \ A+iB\in U(n)\,\right\}
  \cong U(n).
\end{equation*}
Thus we get the biholomorphic map
\begin{equation*}
G/K \lrt \BH_n, \qquad gK \mapsto g\!\cdot\! iI_n,  \quad g\in G.
\end{equation*}
It is well known that $\BH_n$ is an Einstein-K{\"a}hler Hermitian symmetric manifold.
Since $\BH_n$ is K{\"a}hler, it is a symplectic manifold.

\vskip 0.21cm For $Z=(z_{ij})\in\BH_n,$ we write $Z=X+iY$
with $X=(x_{ij}),\ Y=(y_{ij})$ real. We put $dZ=(dz_{ij})$ and
$d{\overline Z}=(d{\overline z}_{ij})$. We
also put
$${ {\!\!\partial}\over {\partial Z} }
=\,\left(\,
{ {1+\delta_{ij}}\over 2}\, { {\!\!\partial}\over {\partial z_{ij} }
} \,\right) \qquad\text{and}\qquad
{ {\!\!\partial}\over {\partial {\overline
{Z}}_{ij} } }
=\,\left(\, {
{1+\delta_{ij}}\over 2}\, { {\!\!\partial}\over {\partial {\overline
{z}}_{ij} } } \,\right).$$ C. L. Siegel \cite{Si1} introduced
the symplectic metric $ds_{n;A}^2$ on $\BH_n$ invariant under the action
\eqref{(2.1)} of $Sp(2n,\BR)$ that is given by
\begin{equation*}
ds_{n;A}^2=A\cdot Tr (Y^{-1}dZ\, Y^{-1}d{\overline Z}),\qquad A>0.
\end{equation*}
It is known that the metric $ds_{n;A}^2$ is a K{\"a}hler-Einstein metric.
H. Maass \cite{M1} proved that its Laplace operator $\Delta_{n;A}$ is given by
\begin{equation*}
\Delta_{n;A}=\,{4\over A}\cdot Tr \left(\,Y\,
{}^{{}^{{}^{{}^\text{\scriptsize $t$}}}}\!\!\!
\left(Y\PZB\right)\PZ\right).
\end{equation*} And
\begin{equation*}
dv_n(Z)=(\det Y)^{-(n+1)}\prod_{1\leq i\leq j\leq n}dx_{ij}\,
\prod_{1\leq i\leq j\leq n}dy_{ij}
\end{equation*}
is a $Sp(2n,\BR)$-invariant volume element on $\BH_n$\,(cf.\,\cite[p.\,130]{Si2}).

\vskip 3mm
Let $\mathbb D(\BH_n)$ be the algebra of all differential operators on $\BH_n$ invariant
under the action \eqref{(2.1)}. Then according to Harish-Chandra
\cite{HC1, HC2},
$$  \mathbb D (\BH_n)= \BC [ H_1,H_2,\cdots,H_n ],$$
where $H_1,H_2,\cdots,H_n$ are algebraically independent invariant differential operators on $\BH_n$.
We note that $n$ is the rank of $\BH_n$, i.e., the rank of $Sp(2n,\BR).$
That is, $\mathbb D (\BH_n)$ is a commutative algebra that is finitely generated by $n$
algebraically independent invariant differential operators on $\BH_n$.
Maass \cite[pp.\,103--121]{M2} found the explicit algebraically independent generators
$H_1,H_2,\cdots,H_n$. Let $\frak g_\BC$ be the complexification of the Lie algebra of $G$.
It is known that $\BD(\BH_n)$ is isomorphic to the center of the universal enveloping algebra of $\fg_{\BC}$\,(cf.\,\cite[Chapter II]{He}).

\vskip 2mm
For each positive integer $n\in \BZ^+,$ we let
\begin{equation}\label{(4.1)}
  G_n:=\,Sp(2n,\BR),\ \ \ K_n:=\,U(n),\ \ \ \G_n:=\,Sp(2n,\BZ)
\end{equation}
be the symplectic group of degree $n$, the unitary group of degree $n$ and the Siegel modular group of degree $n$ respectively. {\it We\ fix\ a\ finite\ dimensional\ complex\ vector\ space} $V$.
And we put $G_0=K_0=\G_0=\{{\rm identity}\}.$
For any two integers $m,n\in\BZ^+$ with $m<n$, we define the
monomorphism
\begin{equation}\label{(4.2)}
u_{m,n}:K_m\lrt K_n
\end{equation}
by
\begin{equation*}
  u_{m,n}(A):=\begin{pmatrix}
               A & 0 \\
               0 & I_{n-m}
             \end{pmatrix},\qquad A\in K_m:=U(m).
\end{equation*}
Let
\begin{equation*}
  K_\infty=\,U(\infty):=\varinjlim_n K_n
\end{equation*}
be the inductive limit of the directed system $(K_n,u_{m,n}).$
Let $\rho_\infty:=\,(\rho_n)$ be a stable representation of $K_\infty$,
that is,
\begin{equation*}
  \rho_\infty:=\varinjlim_n \rho_n,
\end{equation*}
where $\rho_n$ is a rational representation of $K_n$ on $V$ and
for any two positive integers $m,n\in\BZ^+$ with $m<n$,
\begin{equation*}
  \rho_{m}(A):=\rho_n \begin{pmatrix}
               A & 0 \\
               0 & I_{n-m}
             \end{pmatrix},\qquad A\in K_m.
\end{equation*}
For each positive integer $n\in \BZ^+,$ we let ${\bf A}(\rho_n,\Gamma_n)$ be the vector space of
automorphic forms of type $(\rho_n,\G_n).$  See section 3 for the definition of
automorphic forms of type $(\rho_n,\G_n).$
For each positive integer $n\in \BZ^+,$ we
extend $\rho_n$ to the complexification $GL(n,\BC)$ of $K_n$ and also denote by $\rho_n$
the extension of $\rho_n$ to $GL(n,\BC).$
We note that each coset space $G_n/K_n\,(\,n\in \BZ^+\,)$ is an Einstein-K{\"a}hler
Hermitian symmetric space of noncompact type and is biholomorphic to the Siegel upper half plane
\begin{equation*}
\BH_n:=\,\left\{\,Z\in \BC^{(n,n)}\,\vert\ Z=\,^tZ,\ \ \text{Im}\,Z>0\ \right\} \qquad
({\rm see\ Notations}).
\end{equation*}
We recall that $G_n$ acts on $\BH_n$ transitively by
\begin{equation*}
 g\cdot Z:=\,(AZ+B)(CZ+D)^{-1},
\end{equation*}
where $g=\begin{pmatrix} A & B\\ C & D \end{pmatrix}\in G_n$ and $Z\in \BH_n.$ Thus
$G_n/K_n$ is identified with $\BH_n$ via
$$G_n/K_n\ni gK_n\longmapsto g\cdot (iI_n)\in \BH_n.$$
Now for each positive integer $n\in \BZ^+,$ we define the automorphic factor
$J_n:G_n\times \BH_n\lrt GL(V)$ by
\begin{equation}\label{(4.3)}
  J_n(g,Z):=\,\rho_n(CZ+D),
\end{equation}
where $g=\begin{pmatrix} A & B \\ C & D \end{pmatrix} \in G_n$ and $Z\in \BH_n.$

\vskip 3mm
For each positive integer $n\in \BZ^+,$ we denote by $[\G_n,\rho_n]$ the vector space of
Siegel modular forms on $\BH_n$ of type $\rho_n.$ We recall that a Siegel modular form $f$ in $[\G_n,\rho_n]$ is a holomorphic function $f:\BH_n\lrt V$ satisfying the condition
\begin{equation}\label{(4.4)}
f(\g \cdot Z)\,=\,\rho_n(CZ+D)f(Z)\qquad
{\rm for\ all}\ \g=\begin{pmatrix} A & B\\ C & D
\end{pmatrix} \in \G_n\ {\rm and}\ Z\in \BH_n.
\end{equation}
For $n=1,\ f$ requires a cuspidal condition, that is, $f$ is
bounded in any domain $Y\geq Y_0>0$.
For all two positive integers $m,n\in \BZ^+$ with $m< n$, we have the well-known
classical {\sf Siegel\ operator} $\Phi_{m,n}:[\G_n,\rho_n]\lrt
[\G_m,\rho_m]$ defined by
\begin{equation}\label{(4.5)}
\left(\,\Phi_{m,n}f\,\right)(Z):=\,\lim_{t\lrt\infty}\,
f\begin{pmatrix} Z & 0 \\ 0 & it I_{n-m}\end{pmatrix},\quad
Z\in \BH_m.
\end{equation}
We observe that \eqref{(4.5)} is well defined and is a linear mapping.
\vskip 2mm
For an element $F\in {\bf A}(\rho_n,\Gamma_n),$ we define the function $P_nF$ on $\BH_n$ by
\begin{equation}\label{(4.6)}
(P_nF)(g\cdot iI_n):=\,J_n(g,iI_n)F(g),
\end{equation}
where $g\in G_n.$

\begin{lemma}\label{lem:4.1}
If $F\in A(\rho_n,\Gamma_n),$ then $P_nF$ satisfies the condition (4.7).
\end{lemma}
\noindent
{\it Proof.} For any $Z\in \BH_n,$
suppose $Z=g\!\cdot\!iI_n= {\tilde g}\!\cdot\!iI_n, \ g,{\tilde g}\in G_n$ with
$\tilde g =gk,\ k\in K_n.$
We write $k=\begin{pmatrix}
              a & -b \\
              b & \,a
            \end{pmatrix}\in K_n.$
Then for any $F\in A(\rho_n,\Gamma_n),$
\begin{align*}
(P_n F)(Z) &=(P_n F)({\tilde g}\!\cdot\! iI_n) \\
   &= J_n ({\tilde g},iI_n) F(gk) \\
   &= J_n (g,iI_n) J_n (k,iI_n) \rho_n (k)^{-1} F(g) \quad
   ({\rm by\ the\ condition}\ {\rm (AF1)})\\
   &= J_n (g,iI_n) \rho_n (ib+a) \rho_n (ib+a)^{-1} F(g)\\
   &= J_n (g,iI_n) F(g)=(P_n F)(g\!\cdot\! iI_n).
\end{align*}
Thus $P_n F$ is well defined.

\vskip 2mm\noindent
We put $Z=g\!\cdot\!iI_n$ for some $g\in G_n.$
For any $\gamma=\begin{pmatrix}
                  A & B \\
                  C & D
                \end{pmatrix}\in \G_n$
and for any $F\in A(\rho_n,\Gamma_n),$
\begin{align*}
(P_nF)(\gamma\!\cdot\! Z)&=(P_n F) ((\gamma g)\!\cdot\! iI_n)\\
   &=J_n (\g g,iI_n) F(\g g) \\
   &=J_n (\g, g\!\cdot\! iI_n)\, J_n(g,iI_n)F(g) \\
   &=J_n(\gamma,Z) (P_n F)(Z)\\
   &=\rho_n (CZ+D) (P_n F)(Z).
\end{align*}
Therefore $P_nF$ satisfies the condition \eqref{(4.4)}.
\hfill $\square$

\vskip 3mm
For an element $f\in [\G_n,\rho_n],$ we define the function $Q_nf$ on $G_n$ by
\begin{equation}\label{(4.7)}
(Q_nf)(g):=\,J_n(g,iI_n)^{-1}f(g\cdot iI_n),\quad g\in G_n.
\end{equation}

\begin{lemma}\label{lem:4.2}
If $f\in [\G_n,\rho_n],$ then $Q_nf$ is contained in ${\bf A}(\rho_n,\Gamma_n)$.
\end{lemma}
\noindent
{\it Proof.} For any $\g\in \G_n,\ g\in G_n$ and $k\in K_n,$
\begin{align*}
(Q_n f)(\g gk) &= J_n (\g gk,iI_n)^{-1} f((\g gk)\!\cdot\! iI_n)\\
   &= \left( J_n (\g g,k\!\cdot\! iI_n) J_n (k,iI_n)\right)^{-1} f((\g g)\!\cdot\! iI_n) \\
   &= \left( J_n (\g ,g\!\cdot\! iI_n) J_n (g,iI_n) J_n (k,iI_n)\right)^{-1}
            f(\g\!\cdot\!(g\!\cdot\! iI_n)) \\
   &= J_n (k,iI_n)^{-1} J_n (g,iI_n)^{-1}  J_n (\g ,g\!\cdot\! iI_n)^{-1}
      f(\g\!\cdot\!(g\!\cdot\! iI_n))   \\
   &= J_n (k,iI_n)^{-1} J_n (g,iI_n)^{-1}
      f(g\!\cdot\! iI_n)  \\
   &= J_n (k,iI_n)^{-1} (Q_n f)(g) \\
   &= \rho_n (k)^{-1} (Q_n f)(g).
\end{align*}
Thus the condition (AF1) is satisfied. The condition (AF2) is satisfied because $f$,
being of a fixed $K_n$-type $\rho_n$, corresponds to a vector in a Harish-Chandra module,
which is by definition $Z(\mathfrak{g}_\BC)$--finte. 
The condition (AF3) follows from the fact that $f(Z)$ is bounded in
any domain $\left\{ Z\in \BH_n\,\vert\ Z=X+i\,Y,\ Y\geq Y_0>0\,\right\}$
for some positive definite matrix $Y_0$ of degree $n$.
\hfill $\square$

\vskip 1mm
From now on, we denote by
${\bf A}_h(\rho_n,\G_n)$ the image of $[\G_n,\rho_n]$ under $Q_n.$

\vskip 3mm
For all $m,n\in \BZ^+$ with $m<n,$ we define the {\sf Siegel\ operator} $L_{m,n}$
on ${\bf A}_h(\rho_n,\G_n)$ by
\begin{equation}\label{(4.8)}
\left(\,L_{m,n}f\,\right)(g):=\,J_m(g,iI_m)^{-1}\lim_{t\lrt\infty}
J_n(g_t,iI_n)f(g_t),\qquad g\in G_m,
\end{equation}
where $f\in {\bf A}_h(\rho_n,\G_n)$ and $g_t\in G_n$ is defined by
\begin{equation}\label{(4.9)}
g_t:=\,\begin{pmatrix} A & 0 & B & 0\\
0 & t^{1/2}I_{n-m} & 0 & 0 \\
C & 0 & D & 0 \\
0 & 0 & 0 & t^{-1/2}I_{n-m} \end{pmatrix},\quad t>0
\end{equation}
for all $g=\begin{pmatrix} A & B\\ C & D \end{pmatrix} \in G_m.$

\begin{proposition}\label{prop:4.3}
The limit in \eqref{(4.8)} exists and $L_{m,n}$ is a
linear mapping of ${\bf A}_h(\rho_n,\G_n)$ into ${\bf A}_h(\rho_m,\G_m).$
\end{proposition}
\noindent
{\it Proof.}
Let $F=Q_n f\in {\bf A}_h (\rho_n)$ for some $f\in [\G_n,\rho_n].$
Let $g=\begin{pmatrix}
         A & B \\
         C & D
       \end{pmatrix} \in G_m $ and let $g_t \in G_n\,(t>0)$
be the element in $G_n$ given by the formula \eqref{(4.9)}. Then we have
\begin{align*}
(L_{m,n}F)(g) &= J_m(g,iI_m)^{-1} \lim_{t\lrt\infty}
J_n(g_t, iI_n) F(g_t) \\
   &= J_m(g,iI_m)^{-1} \lim_{t\lrt\infty}
J_n(g_t, iI_n) (Q_n f)(g_t)   \\
   &= J_m(g,iI_m)^{-1} \lim_{t\lrt\infty}
J_n(g_t, iI_n) J_n(g_t, iI_n)^{-1} f(g_t\!\cdot\! iI_n) \quad (by\ {\rm (4.7)})\\
   &= J_m(g,iI_m)^{-1} \lim_{t\lrt\infty}
      f \begin{pmatrix}
          (Ai+B)(Ci+D)^{-1} & 0 \\
           0 & i\,tI_{n-m}
        \end{pmatrix}.
\end{align*}
Since $f$ is an element of $[\G_n,\rho_n]$,
the limit
$$\lim_{t\lrt\infty}
      f \begin{pmatrix}
          (Ai+B)(Ci+D)^{-1} & 0 \\
           0 & i\,tI_{n-m}
        \end{pmatrix}=(\Phi_{m,n}f)(g\cdot iI_m)$$
exists and is an element of $[\G_m,\rho_m]$. Here $\Phi_{m,n}$ is the Siegel operator
defined by the formula \eqref{(4.4)}.
Thus the limit in \eqref{(4.8)} exists. On the other hand,
\begin{align*}
  (L_{m,n}F)(g) &= J_m (g,iI_m)^{-1}(\Phi_{m,n}f)(g\!\cdot\! iI_m) \\
   &= \left( Q_m (\Phi_{m,n}f)\right)(g).
\end{align*}
Therefore $L_{m,n}F$ is an element of ${\bf A}_h(\rho_m).$
Hence $L_{m,n}$ is a
linear mapping of ${\bf A}_h(\rho_n)$ into ${\bf A}_h(\rho_m).$
\hfill $\square$

\vskip 5mm
Let ${\frak g}_n$ be the Lie algebra of $G_n$ and ${\frak g}_n^{\BC}$ its complexification. Then
\begin{equation*}
{\frak g}_n^{\BC}=\left\{\begin{pmatrix} A & B\\ C & -^t\!A\end{pmatrix}
\in \BC^{(2n,2n)}\,
\bigg|\ B=\,^tB,\ \ C=\,^tC\ \right\}.
\end{equation*}
We let $\widehat{J}_n:=iJ_n$ with
$J_n=\begin{pmatrix} 0 & I_n\\ -I_n & 0 \end{pmatrix}.$
We define an involution $\s_n$ of $G_n$ by
$$\s_n(g):=\widehat{J}_ng\widehat{J}_n^{-1},\ \ \ g\in G_n.$$
The differential map $d\s_n={\rm Ad}\,(\widehat{J}_n)$ of $\s_n$ extends complex linearly
to the complexification ${\frak g}_n^{\BC}$ of ${\frak g}_n.
\ \text{Ad}\,(\widehat{J}_n)$ has
1 and -1 as eigenvalues. The $(+1)$-eigenspace of $\text{Ad}\,(\widehat{J}_n)$
is given by
\begin{equation*}
{\frak k}_n^{\BC}:=\left\{
\begin{pmatrix} A&-B\\ B & A\end{pmatrix}\in \BC^{(2n,2n)}\ \bigg|\
^t\!A+A=0,\ B=\,^tB\ \right\}.
\end{equation*}
We note that ${\frak k}_n^{\BC}$ is the complexification of the Lie algebra
${\frak k}_n$ of a maximal compact subgroup $K_n=G_n\cap SO(2n,\BR)\cong U(n)$ of
$G_n$. The $(-1)$-eigenspace of $\text{Ad}\,(\widehat{J}_n)$ is given by
\begin{equation*}
{\frak p}_n^{\BC}=\left\{ \begin{pmatrix} A & B\\ B & -A
\end{pmatrix}\in \BC^{(2n,2n)}\
\bigg|\ A=\,^t\!A,\ B=\,^tB\ \right\}.
\end{equation*}
We observe that ${\frak p}_n^{\BC}$ is not a Lie algebra.
But ${\frak p}_n^{\BC}$ has the following decomposition
$${\frak p}_n^{\BC}={\frak p}_{n,+}\oplus {\frak p}_{n,-},$$
where
$${\frak p}_{n,+}=\left\{\begin{pmatrix} X & iX \\ iX & -X
\end{pmatrix} \in \BC^{(2n,2n)}\
\bigg|\ X=\,^tX\ \right\}$$
and
$${\frak p}_{n,-}=\left\{\begin{pmatrix} Y & -iY\\ -iY & -Y
\end{pmatrix}\in \BC^{(2n,2n)}\
\bigg|\
Y=\,^tY\ \right\}.$$
We observe that ${\frak p}_{n,+}$ and ${\frak p}_{n,-}$ are abelian subalgebras of
${\frak g}_n^{\BC}.$

\begin{proposition}\label{prop:4.4}
A function $F$ in ${\rm A}_h(\rho_n,\G_n)$ is characterized
by the following properties ${\rm (Sp1)}\!-\!{\rm (Sp3)}$\,:
\begin{equation*}
\ \ {\rm (Sp1)}\qquad\quad  F(\g gk)=\rho_n(k)^{-1}F(g)\quad {\rm for\ all}\ \g\in \G_n,
\ g\in G_n\ {\rm and}\ k\in K_n. \hskip 5cm
\end{equation*}
\begin{equation*}
\ \ {\rm (Sp2)}\qquad\qquad \qquad
X^{-}F=0 \quad {\rm for \ all}\ X^{-}\in {\frak p}_{n,-}.\hskip 8cm
\end{equation*}
\ \ ${\rm (Sp3)}$ \ For any $M\in G_n,$ the function $\psi:G_n\lrt V$ defined by
\begin{equation*}
\psi(g):=\,\rho_n(Y^{-{\frac 12}})F(Mg),\quad \ g\in G_n,\ \ g\cdot iI_n:=X+i\,Y
\end{equation*}
\indent\indent\ \ \ \ \
is bounded in the domain $Y\geq Y_0>0$ for some $Y_0=\,^tY_0>0.$
\end{proposition}
\noindent
{\it Proof.} (Sp1) follows from (AF1) and (Sp2) follows from the fact
that $f$ is holomorphic on $\BH_n$. Since $F\in A_h(\rho_n):={\rm Im}\, Q_n,\
F=Q_nf$ for some $f\in [\G_n,\rho_n].$ Let $Z=X+i\,Y=g\!\cdot\! iI_n\in\BH_n$
with $g\in G_n$. Then for any $M\in G_n,$ we have
\begin{align*}
\psi (g) &= \rho_n(Y^{-{\frac 12}})F(Mg)  \\
   &= \rho_n(Y^{-{\frac 12}})J_n(Mg,iI_n)^{-1} f((Mg)\!\cdot\!(iI_n)) \\
   &= \rho_n(Y^{-{\frac 12}})J_n(g,iI_n)^{-1} J_n(M,Z)^{-1}f(M\!\cdot\!Z).
\end{align*}
If $Y_0$ is sufficiently large, $\psi(g)$ is bounded.
\hfill $\square$

\begin{proposition}\label{prop:4.5}
The mapping $P_n$ and $Q_n$ are compatible with
the Siegel operators $L_{m,n}$ and $\Phi_{m,n}\,(\,m<n\,)$. That is,
for any $m,n\in \BZ^+$ with $m<n$, we have
\begin{equation}\label{(4.10)}
L_{m,n}\circ Q_n\,=\,Q_m\circ \Phi_{m,n}\quad \ \text{on}\ \
[\G_n,\rho_n]
\end{equation}
and
\begin{equation}\label{(4.11)}
P_m\circ L_{m,n}\,=\,\Phi_{m,n}\circ P_n\quad \ \text{on}\ \
A_h(\rho_n,\G_n).
\end{equation}
\end{proposition}
\noindent
{\it Proof.} Let $f\in [\G_n,\rho_n]$ and
$g=\begin{pmatrix}
   A & B \\
   C & D
 \end{pmatrix}\in G_m.$ Let $g_t\in G_n\,(t>0)$ be the matrix defined by
the formula \eqref{(4.9)}. Then
\begin{align*}
  \left( L_{m,n} (Q_n f)\right)(g) &=
  J_m (g,iI_m)^{-1} \lim_{t\lrt \infty} J_n (g_t,iI_n) (Q_n f)(g_t) \\
  &= J_m (g,iI_m)^{-1} \lim_{t\lrt \infty} J_n (g_t,iI_n)
     J_n (g_t,iI_n)^{-1} f(g_t\!\cdot\! iI_n) \\
  &= J_m (g,iI_m)^{-1} \lim_{t\lrt \infty}
     f \begin{pmatrix}
         (Ai+B)(Ci+D)^{-1} & 0 \\
         0 & i\,t\, I_{n-m}
       \end{pmatrix}.
\end{align*}
\noindent
On the other hand,
\begin{align*}
  \left( Q_m(\Phi_{m,n}(f))\right)(g)
  &= J_m (g,iI_m)^{-1} (\Phi_{m,n}(f))(g\!\cdot\! iI_m)\\
  &= J_m (g,iI_m)^{-1} \lim_{t\lrt \infty}
     f \begin{pmatrix}
        g\!\cdot\! iI_m  & 0 \\
         0 & i\,t\, I_{n-m}
       \end{pmatrix}   \\
  &= J_m (g,iI_m)^{-1} \lim_{t\lrt \infty}
     f \begin{pmatrix}
         (Ai+B)(Ci+D)^{-1} & 0 \\
         0 & i\,t\, I_{n-m}
       \end{pmatrix}.
\end{align*}
\noindent
This proves the formula \eqref{(4.10)}.
\vskip 2mm
Let $F\in A_h (\rho_m,\G_m)$ and let
$Z=g\cdot iI_m \in \BH_m$ with $g\in G_m.$
Let $g_t\in G_n\,(t>0)$ be the matrix defined by the formula \eqref{(4.9)}. Then
\begin{align*}
  \left( P_m(L_{m,n}(F))\right)(Z)
  &= J_m (g,iI_m) (L_{m,n}(F))(g)\\
  &= J_m (g,iI_m) J_m (g,iI_m)^{-1}
  \lim_{t\lrt \infty}
     J_n (g_t,iI_n) F(g_t) \\
  &=\lim_{t\lrt \infty}
     J_n (g_t,iI_n) F(g_t).
\end{align*}
On the other hand,
\begin{align*}
  \left( \Phi_{m,n}(P_n F)\right)(Z)
  &= \lim_{t\lrt \infty}
     (P_n F) \begin{pmatrix}
        Z  & 0 \\
         0 & i\,t\, I_{n-m}
       \end{pmatrix}   \\
  &=  \lim_{t\lrt \infty} (P_n F)(g_t\!\cdot \! iI_n)\\
  &=  \lim_{t\lrt \infty} J_n (g_t,iI_n) F(g_t).
\end{align*}
Therefore the formula \eqref{(4.11)} is proved. \hfill $\square$

\vskip 2mm
We set
\begin{equation*}
\G_\infty=\varinjlim_n \G_n=\varinjlim_n Sp(2n,\BZ)\qquad {\rm and}
\qquad \rho_\infty:=\,\varinjlim_n  \rho_n.
\end{equation*}
Using the Siegel operator $\Phi_{m,n},$ we define the inverse limit
\begin{equation}\label{(4.12)}
[\G_\infty,\,\rho_\infty]:=\,\varprojlim_n  [\G_n,\,\rho_n]
\end{equation}

For $n\in \BZ^+,$ we put
$$\mathbb M_n:=\,\bigoplus_{\tau}[\G_n,\,\tau],$$
where $\tau$ runs over all isomorphism classes of irreducible rational
finite dimensional representations of the general linear group $GL(n,\BC)$ of
degree $n$. For $n=0,$ we set $\mathbb M_0:=\BC.$ For an irreducible finite
dimensional representation $\tau=(\la_1,\la_2,\cdots,\la_n)$
of $GL(n,\BC)$ with $\la_1\geq \la_2\cdots\geq \la_n,\ \la_i\in \BZ\,(\,
1\leq i\leq n\,),$ the integer $k(\tau):=\la_n$ is called the
{\sf weight} of $\tau.$ Here $\tau=(\la_1,\la_2,\cdots,\la_n)$ denotes
the irreducible finite dimensional representation of $GL(n,\BC)$ with
highest weight $(\la_1,\la_2,\cdots,\la_n)\in\BZ^n.$

\vskip 2mm
For a positive integer $n\in \BZ^+,$ we define
$$\mathbb M_n^{\ast}:=\,\bigoplus_{\tau:\text{even}}[\G_n,\tau],$$
where $\tau$ runs over all isomorphism classes of irreducible rational
finite dimensional {\sf even} representations of $GL(n,\BC)$ such that the
highest weight $\lambda(\tau)$ of $\tau$ is even, i.e., $\lambda(\tau)
\in (2\BZ)^n$. For $n=0,$ we also set
$\mathbb M_0^{\ast}:=\,\BC.$ Clearly for any two positive integers $m,n$ with $m<n$,
the Siegel operator $\Phi_{m,n}$ maps
$\mathbb M_n$\,(\,resp.\,$\mathbb M_n^{\ast}$\,) into $\mathbb M_m$\,(\,resp.\,$\mathbb M_m^{\ast}$\,).

\vskip 2mm
We let
\begin{equation*}
\mathbb M:=\,\lim_{\begin{subarray}{c} \longleftarrow\\ ^n \end{subarray}} \mathbb M_n\qquad \text{and}\qquad
\mathbb M^{\ast}:=\,\lim_{\begin{subarray}{c} \longleftarrow\\ ^n \end{subarray}} \mathbb M_n^{\ast}.
\end{equation*}
It is easy to see that both $\mathbb M$ and $\mathbb M^{\ast}$ have
commutative ring structures compatible with the Siegel operators
$\Phi_{\ast,\ast}.$ Obviously $\mathbb M^{\ast}$ is a subring of
$\mathbb M.$

\vskip 3mm
Now we obtained the following result.
\begin{proposition}\label{prop:4.6}
$$\mathbb M\,=\,\bigoplus_{\rho_\infty}[\G_\infty,\,\rho_\infty],$$
where $\rho_\infty$ runs over all stable representations.
\end{proposition}

\vskip 3mm
\begin{definition}\label{def:4.7}
Elements of $\mathbb M$ are called
{\sf stable\ modular\ forms} and elements of $\mathbb M^{\ast}$ are
called {\sf even\ stable\ modular\ forms.}
\end{definition}

\begin{remark}\label{rk:4.8}
As mentioned before in the introduction, the concept of stable modular forms was first introduced by
E. Freitag\,\cite{Fr2}. Thereafter the study of stable modular forms was intensively investigated by R. Weissauer\,\cite{We}.
\end{remark}

\vskip 3mm
Now we give an example of stable modular forms.
\begin{definition}\label{def:4.9}
A pair $(\Lambda,Q)$ is called a {\sf quadratic form} if $\Lambda$ is
a lattice and $Q$ is an integer-valued bilinear symmetric form on
$\Lambda$. The {\sf rank} of $(\Lambda,Q)$ is defined to be the rank of
$\Lambda$. For $v\in \Lambda$, the integer $Q(v,v)$ is called the {\sf norm}
of $v$. A quadratic form $(\Lambda,Q)$ is said to be {\sf even} if
$Q(v,v)$ is even for all $v\in \Lambda$.
A quadratic form $(\Lambda,Q)$ is said to be {\sf unimodular} if $\det (Q)=1.$
\end{definition}

\begin{definition}\label{def:4.10}
Let $(\Lambda, Q)$ be an even unimodular positive definite quadratic form of rank $r$.
For a positive integer $n$, the theta series $\theta_{Q,n}$ associated to $(\Lambda, Q)$ is defined to be
\begin{equation*}
\theta_{Q,n}(\tau):=\sum_{x_1,\cdots,x_n\in \Lambda}\exp \left( \pi i \sum_{p,q=1}^n Q(x_p,x_q)\tau_{pq}\right),\qquad \tau=(\tau_{pq})\in \BH_n.
\end{equation*}
\end{definition}

\begin{proposition}\label{prop:4.11}
Let $(\Lambda, Q)$ be an even unimodular positive definite quadratic form of rank $r$. Then the collection of all theta series associated to $(\Lambda, Q)$
\begin{equation*}
\Theta_Q:=\left( \theta_{Q,n} \right)_{n\geq 0}
\end{equation*}
is a stable modular form of weight ${\frac r2}.$
\end{proposition}
\noindent
{\it Proof.} We note that $r\equiv 0\,({\rm mod}\,8)$ (cf.\,\cite{Ser}). It is well known that $\theta_{Q,n}(\tau)$ is a Siegel modular form on $\BH_n$ of weight ${\frac r2}$ (cf.\,\cite{Fr3}). We easily see that
\begin{equation*}
\Phi_{m,n}(\theta_{Q,n})=\theta_{Q,m}\qquad {\rm for\ all}\ m,n\  {\rm with}\ m<n.
\end{equation*}
Therefore the collection
$\Theta_Q=\left( \theta_{Q,n} \right)_{n\geq 0}$
is a stable modular form of weight ${\frac r2}.$
\hfill $\square$


\vskip 7mm
\noindent
{\bf Example B.\ Stable automorphic forms for $SL(\infty,\BR)$}
\vskip 3mm
First of all, we provide some geometric properties on the symmetric manifold
$SL(n,\BR)/SO(n,\BR)$ that is important geometrically and number theoretically.
\vskip 2mm
\begin{definition}\label{def:4.12}
For any positive integer $n\geq 2$. we define ${\mathfrak H}_n$ to be the set of all
$n\times n$ real matrices of the form $z=x\cdot y$, where
$$x=\begin{pmatrix}
      1 & x_{12} & x_{13} &\cdots & x_{1n} \\
      0 & 1 &  x_{22} & \cdots & x_{2n} \\
      0 & 0 & \ddots & \vdots & \vdots  \\
      {0} & {0} & {0} & {0} & 1
    \end{pmatrix}$$
and
$$y={\rm {diag}}(y_1y_2\cdots y_{n-1},y_1y_2\cdots y_{n-2},\cdots,y_1,1)$$
with $x_{ij}\in \BR$ for $1\leq i < j \leq n$ and $y_i > 0 $ for $1\leq i\leq n-1.$
\end{definition}
Let
$$\mathbb{X}_n:=\{\, Y=\,{}^tY >0,\ \det\,(Y)=1\,\}.$$
We can show that ${\mathfrak H}_n$ is diffeomorphic to ${\mathbb X}_n$.
In fact, we have the Iwasawa decomposition
\begin{equation*}
  GL(n,\BR)={\mathfrak H}_n\cdot O(n)\cdot Z_n,
\end{equation*}
where $Z_n(\cong \BR^{\times})$ is the center of $GL(n,\BR)$
\,(cf.\,\cite[Proposition 1.2.6, pp.\,11--12]{Go}).
Here
\begin{equation*}
  O(n):=O(n,\BR)=\{ k\in GL(n,\BR)|\ {}^tk k=k\,{}^tk=I_n\,\}
\end{equation*}
denotes the real orthogonal group of degree $n$.
We see easily that
\begin{equation*}
  {\mathfrak H}_n \cong GL(n,\BR)/(O(n)\cdot \BR^{\times})\cong
  SL(n,\BR)/SO(n,\BR)\cong {\mathbb X}_n,
\end{equation*}
where $\cong$ denotes the diffeomorphism.
\vskip 2mm
It is seen that $GL(n,\BR)$ acts on ${\mathfrak H}_n$ by left translation
\,(cf.\,\cite[Proposition 1.2.10, p.\,14]{Go}).
Then we obtain
\begin{equation*}
 {\mathfrak X}_n:=SL(n,\BZ)\backslash SL(n,\BR)/SO(n) \cong SL(n,\BZ) \backslash
  GL(n,\BR)/(O(n)\cdot \BR^{\times}),
\end{equation*}
where $SO(n):=SO(n,\BR)=SL(n,\BR)\cap O(n).$
Therefore we obtain the following isomorphism
\begin{equation*}
  {\mathfrak X}_n\cong SL(n,\BZ)\backslash {\mathfrak H}_n.
\end{equation*}

\begin{proposition}\label{prop:4.13}
Let $n\geq 2.$ Following the coordinates of Definition 4.4, we put
\begin{equation*}
  d^*x=\prod_{1\leq i<j\leq n} dx_{ij}\quad {\rm and}\quad
  dy^*=\prod_{k=1}^{n-1} y_k^{-n(n-k)-1}dy_k.
\end{equation*}
Then
\begin{equation*}
  d^*z=d^*x\cdot d^*y
\end{equation*}
is the left $SL(n,\BR)$-invariant volume element on
${\mathbb X}_n\cong {\mathfrak H}_n.$
\end{proposition}
\begin{proof}
The proof can be found in \cite[Proposition 1.5.3, pp.\,25--26]{Go}.
\end{proof}

\vskip 2mm
Following earlier work of Minkowski, Siegel\,\cite{Si} calculated the volume
${\bf Vol}(\G^n\ba {\mathfrak H}_n)$ which is expressed in terms of
$$
\zeta(2)\cdot \zeta(3)\cdots \zeta(n).
$$
\begin{theorem}\label{thm:4.14}
Let $n\geq 2.$ Then the volume ${\bf Vol}(\G^n\ba {\mathfrak H}_n)$ of
$\G^n\ba {\mathfrak H}_n$ is given by
\begin{equation*}
  {\bf Vol}(\G^n\ba {\mathfrak H}_n)=\int_{\G^n\ba {\mathfrak H}_n} d^*z
  =n\cdot 2^{n-1}\cdot \prod_{k=2}^n { {\zeta(k)}\over {{\bf Vol}(S^{k-1})} },
\end{equation*}
where $\G^n=SL(n,\BZ)$ and
\begin{equation*}
  {\bf Vol}(S^{k-1})={ {2(\sqrt{\pi})^k}\over {\Gamma(k/2)} }
\end{equation*}
denotes the volume of the $(k-1)$-dimensional sphere $S^{k-1}$,
$\zeta(k)=\sum_{n=1}^{\infty}{1\over {n^k}}$ is the Riemann zeta function and
$\G(p)$ denotes the Gamma function.
\end{theorem}
\begin{proof}
The proof can be found in \cite[Theorem 1.6.1, pp.\,27--37]{Go} or \cite{Si}.
\end{proof}

\newcommand\POB{ {{\partial}\over {\partial{\overline \Omega}}} }
\newcommand\PX{ {{\partial}\over{\partial X}} }
\newcommand\PY{ {{\partial}\over {\partial Y}} }
\newcommand\PU{ {{\partial}\over{\partial U}} }
\newcommand\PV{ {{\partial}\over{\partial V}} }
\newcommand\PO{ {{\partial}\over{\partial \Omega}} }
\newcommand\PW{ {{\partial}\over{\partial W}} }
\newcommand\PWB{ {{\partial}\over {\partial{\overline W}}} }
\newcommand\OVW{\overline W}
\newcommand\Rg{{\mathfrak R}_n}

\vskip 2mm
For any positive integer $n\geq 1$, we let
\begin{equation*}
  {\mathscr P}_n:=\{Y\in \BR^{(n,n)}\,|\ Y=\,{}^tY>0\,\}
\end{equation*}
be the open convex cone in the Euclidean space $\BR^N$ with $N=\frac{n(n+1)}{2}.$
Then $GL(n,\BR)$ acts ${\mathscr P}_n$ transitively by
\begin{equation}\label{(4.13)}
  g\cdot Y=gY\,^t\!g\qquad {\rm for\ all}\ g\in GL(n,\BR)\ {\rm and}\ Y\in {\mathscr P}_n.
\end{equation}
Since $O(n)$ is the isotopic subgroup of $GL(n,\BR)$ at $I_n$, the symmetric
space $GL(n,\BR)/O(n)$ is diffeomorphoc to ${\mathscr P}_n$.
\vskip 0.3cm
For $Y=(y_{ij})\in {\mathscr P}_n,$ we put
\begin{equation*}
dY=(dy_{ij})\qquad\text{and}\qquad \PY\,=\,\left(\, {
{1+\delta_{ij}}\over 2}\, { {\partial}\over {\partial y_{ij} } }
\,\right).
\end{equation*}

\vskip 0.2cm For a fixed element $A\in GL(n,\BR)$, we put
$$Y_*=A\circ Y=AY\,^t\!A,\quad Y\in {\mathscr P}_n.$$
Then
\begin{equation}\label{(4.14)}
dY_*=A\,dY\,^t\!A \quad \textrm{and}\quad {{\partial}\over {\partial
Y_*}}=\,^t\!A^{-1} \Yd\, A^{-1}.
\end{equation}

\vskip 5mm
We can see easily that
\begin{equation*}  ds^2=\,Tr( (Y^{-1}dY)^2)  \end{equation*}
is a $GL(n,\BR)$-invariant Riemannian metric on ${\mathscr P}_n$ and its
Laplacian is given by
\begin{equation*}
\Delta=   Tr\left( \left( Y\PY\right)^2\right),
\end{equation*}
\noindent where $ Tr(M)$ denotes the trace of a square
matrix $M$. We also can see that
\begin{equation*}
d\mu_n(Y)=(\det Y)^{-{ {n+1}\over2 } }\prod_{i\leq j}dy_{ij}
\end{equation*}
is a $GL(n,\BR)$-invariant volume element on ${\mathscr P}_n$.

\begin{theorem}\label{thm:4.15}
A geodesic $\alpha (t)$ joining $I_n$ and $Y\in {\mathscr P}_n$ has the form
\begin{equation*}
\alpha (t)=\exp (t A[V]),\qquad t\in [0,1],
\end{equation*}
where
\begin{equation*}
Y=(\exp A)[V]=\exp (A[V])=\exp (\,^tVAV)
\end{equation*}
is the spectral decomposition of $Y$, where $V\in O(n,\BR),\ A={\rm diag} (a_1,\cdots,a_n)$ with all $a_j\in \BR.$
The distance of $\alpha (t) \ (0\leq t\leq 1)$ between $I_n$ and $Y$ is
\begin{equation*}
 \left( \sum_{j=1}^{n} a_j^2 \right)^{\frac{1}{2}}.
\end{equation*}
\end{theorem}
\begin{proof} The proof can be found in \cite[pp.\,16-17]{T}.\end{proof}

\vskip 2mm
We consider the following {\sf Maass}-{\sf Selberg\ (differential)\ operators}
$\delta_1,\delta_2,\cdots,\delta_n$ on ${\mathscr P}_n$ defined by
\begin{equation}\label{(4.15)}
\delta_k= Tr\left( \left( Y\Yd \right)^k\right),\quad
k=1,2,\cdots,n,
\end{equation}
By Formula \eqref{(4.14)}, we get
\begin{equation*}
\left( Y_* {{\partial}\over {\partial Y_*}}\right)^i=\,A\,\left(
Y\Yd\right)^i A^{-1}
\end{equation*}

\noindent for any $A\in GL(n,\BR)$. So each $\delta_i\ (1\leq i \leq n)$ is invariant
under the action \eqref{(4.13)} of $GL(n,\BR)$.

\vskip 0.2cm Selberg \cite{Sel} proved the following.

\begin{theorem}\label{thm:4.16}
The algebra ${\mathbb D}({\mathscr P}_n)$ of all $GL(n,\BR)$-invariant differential operators on
${\mathscr P}_n$ is generated by $\delta_1,\delta_2,\cdots,\delta_n.$ Furthermore
$\delta_1,\delta_2,\cdots,\delta_n$ are algebraically independent and
${\mathbb D}({\mathscr P}_n)$ is isomorphic to the commutative ring $\BC[x_1,x_2,\cdots,x_n]$
with $n$ indeterminates $x_1,x_2,\cdots,x_n.$
\end{theorem}

\begin{proof} The proof can be found in \cite[pp.\,64-66]{M2}.\end{proof}

\vskip 2mm
Using the Maass-Selberg operators, Brennecken, Ciardo and Hilgert \cite{BCH} found the
{\it explicit\ generators} $E_1,E_2,\cdots,E_{n-1}$ of the algebra of all
$SL(n,\BR)$-invariant differential operators on ${\mathbb X}_n.$ We observe that
$E_1,E_2,\cdots,E_{n-1}$ are algebraically independent.

\vskip 3mm
For any $\nu=(\nu_1,\nu_2,\cdots,\nu_{n-1}),$ we define the function
$I_\nu:{\mathfrak H}_n\lrt \BC$ by
\begin{equation*}
  I_\nu (z):=\prod_{i=1}^{n-1} \prod_{i=1}^{n-1}y_i^{b_{ij}\nu_j},
\end{equation*}
where
$$  b_{ij}:=\begin{cases}
            \ \ \ ij, & {\rm if }\ \, i+j\leq n \\
            (n-i)(n-j), & {\rm if}\ \, i+j\geq n.
          \end{cases}
$$
We denote by $\BD(\mathfrak{H}_n)$ the algebra of all $SL(n,\BR)$-invariant
differential operators on $\mathfrak{H}_n.$
Then we see that $I_\nu(z)$ is an eigenfunction of ${\mathbb D}({\mathfrak H}_n).$
Let us write
\begin{equation}
  DI_\nu(z)=\lambda_D\cdot I_\nu(z)\quad {\rm for\ every}\ D\in {\mathbb D}({\mathfrak H}_n).
\end{equation}
\noindent
Since
$$
\lambda_{D_1D_2}=\lambda_{D_1}\lambda_{D_2}\qquad {\rm for\ all}\ D_1,D_2\in
{\mathbb D}({\mathfrak H}_n).
$$
The function $\lambda_D$ (viewed as a function of $D$) is a character of
${\mathbb D}({\mathfrak H}_n$ which is called the {\sf Harish}-{\sf Chandra\ character}.

\vskip 3mm
Following Goldfeld\,(cf.\,\cite[Definition 5.1.3, pp.\,115--116]{Go}, the notion of
a Maass form is defined in the following way.
\begin{definition}\label{def:4.17}
Let $n\geq 2$ and $\Gamma^n=SL(n,\BZ).$ For any $\nu=(\nu_1,\nu_2,\cdots,\nu_{n-1})\in \BC^{n-1},$ a smooth
$f:\G^n\ba {\mathfrak H}_n\lrt \BC$ is said to be a {\sf Maass\ form} for $\G^n$ of
type $\nu$ if satisfies the following conditions (M1)--(M3)\,:
\vskip 2mm\noindent
\ \ ${\rm (M1)}\ \ F(\gamma z)=f(z)\quad {\rm for\ all}\ \gamma\in \G^n,\ z\in {\mathfrak H}_n.$
\vskip 2mm\noindent
\ \ ${\rm (M2)}\ \ Df(z)=\lambda_D f(z)\quad {\rm for\ all}\ D\in {\mathbb D}({\mathfrak H}_n)$
\ {\rm with}\ {\rm eigenvalue}\ $\lambda_D$\ {\rm given\ by}\ (4.16).
\vskip 2mm\noindent
\ \ ${\rm (M3)}\ \ \int_{\G^n\cap U\ba U} f(uz)du=0$
\vskip 1mm
\indent
\ \ \ \ \ \ \ \ for all upper triangular groups $U$ of the form
$$U=\left\{
\begin{pmatrix}
  I_{r_1} & * & * & * \\
  0 & I_{r_2} & * & * \\
  0 & 0 & \ddots & * \\
  0 & 0 & 0 & I_{r_b}
\end{pmatrix}
\right\}
$$
\ \ \ \ \ with $r_1+r_2+\cdots+r_b=n.$ Here $I_r$ denotes the $r\times r$ identity matrix and
$*$ denotes
\par \ \ \ \ \ arbitrary real matrices.
\end{definition}
\begin{remark}\label{rk:4.18}
In \cite{Go}, Dorian Goldfeld studied
Whittaker functions associated with Maass forms, Hecke operators for $\G^n$,
the Godement-Jacquet $L$-function for $\G^n$, Eisenstein series for $\G^n$ and
Poincar{\'e} series for $\G^n$.
\end{remark}

\vskip 2mm
Let
\begin{equation*}
  G^n=SL(n,\BR),\quad K^n=SO(n)\quad {\rm and}\quad \G^n=SL(n,\BZ).
\end{equation*}
Let
\begin{equation*}
  {\mathbb X}_n:=\left\{ Y\in \BR^{(n,n)}\,|\ Y=\,{}^tY>0,\ \ \det Y=1\,
  \right\}
\end{equation*}
be a symmetric space associated to $G^n$. Indeed, $G^n$ acts on ${\mathbb X}_n$ transitively by
\begin{equation}\label{(4.17)}
  g\circ Y= gY\,{}^tg\qquad {\rm for\ all}\ g\in G^n\ {\rm and}\ Y\in {\mathbb X}_n.
\end{equation}
Thus ${\mathbb X}_n$ is a smooth manifold diffeomorphic to the symmetric space
$G^n/K^n$ through the bijective map
\begin{equation*}
  G^n/K^n \lrt {\mathbb X}_n,\qquad gK^n \longmapsto g\circ I_n=g\,{}^tg
  \quad {\rm for\ all}\ g\in G^n.
\end{equation*}

An {\sf automorphic form} for $\G^n$ is defined to be a real analytic function
$f\!:\!{\mathbb X}_n\lrt \BC$ satisfying the following properties (SL1)--(SL3)\,:
\vskip 2mm
(SL1) $f$ is an eigenfunction for all
$G^n$-invariant\ differential\ operators on ${\mathbb X}_n$.
\vskip 2mm
(SL2) \ \ $f(\g Y\,{}^t\g)=f(Y)\quad {\rm for\ all}\ \g\in\G^n
\ {\rm and}\ Y\in {\mathbb X}_n.$
\vskip 2mm
(SL3)\ There\ exist\ a constant  $C>0$ and $s\in \BC^{n-1}$
with  $s=(s_1,\cdots,s_{n-1})$ such that  \\
\indent \ \ \ \ \ \ \ \ \
$|f(Y)| \leq C\,|p_{-s}(Y)|$ as the upper left determinants
$\det Y_j\lrt \infty\ (1\leq j\leq n-1),$ \\
\indent \ \ \ \ \ \ \ \ \ where
\begin{equation*}
  p_{-s}(Y):=\prod_{j=1}^{n-1} (\det Y_j)^{-s_j}
\end{equation*}
\indent \ \ \ \ \ \ \ \ \
is the Selberg's power function\,(cf.\,\cite{Sel, T}).

\vskip 3mm
We denote by ${\bf A}(\G^n)$ the space of all automorphic forms for $\G^n.$
A {\sf cusp form} $f\in {\bf A}(\G^n)$ is defined to be an automorphic form
for $\G^n$ satisfying the following conditions\,:
\begin{equation*}
\int_{X\in (\BR/\BZ)^{(j,n-j)}}
f \left( Y\left[ \begin{pmatrix}
                   I_j & X \\
                   0 & I_{n-j}
                 \end{pmatrix}\right]\right)dX=0,
                 \quad 1\leq j\leq n-1.
\end{equation*}
We denote by ${\bf A}_0(\G^n)$ the space of all cusp forms for $\G^n.$

K
\begin{definition}\label{def:4.19}
Let $f\in {\bf A}(\G^n)$ be an automorphic form for $\G^n$ with eigenvalues
determined by $s=(s_1,\cdots,s_{n-1})\in \BC^{(n-1)}$.
We set
\begin{equation*}
  \xi_1={\frac{1}{n-1}} \sum_{n=2}^{n-1} (n-k)s_k.
\end{equation*}
We define, for any $f\in {\bf A}(\G^n)$,
\begin{equation}\label{(4.18)}
  \mathfrak L_n f (W):=\lim_{v\lrt \infty} v^{-s_1-\xi_1}f(Y),
  \quad v\in \BR,\ W\in {\mathbb X}_{n-1},\ Y\in {\mathbb X}_n,
\end{equation}
where $Y,\ v,\ W$ are determined by the unique decomposition of $Y$
given by
\begin{equation}\label{(4.19)}
  Y=\begin{pmatrix}
      1 & 0 \\
      x & I_{n-1}
    \end{pmatrix}
    \begin{pmatrix}
      v^{-1} & 0 \\
      0 & v^{\frac{1}{n-1}}W
    \end{pmatrix}
    \begin{pmatrix}
      1 & {}^tx \\
      0 & I_{n-1}
    \end{pmatrix},\quad x\in \BR^{(n-1,1)}.
\end{equation}
\end{definition}
D. Grenier \cite{Gr} defined the formula \eqref{(4.18)} and proved the following result.
\begin{theorem}\label{thm:4.20}
If $f\in {\bf A}(\G^n)$, then $\mathfrak L_n f\in {\bf A}(\G^{n-1}).$ Thus
$\mathfrak L_n$ is a linear mapping of ${\bf A}(\G^n)$ into ${\bf A}(\G^{n-1})$.
Moreover if $f\in {\bf A}_0(\G^n)$ is a cusp form, then $\mathfrak L_n f=0.$
In general, $\ker \mathfrak L_n\neq {\bf A}_0(\G^n).$
\end{theorem}
\noindent
{\it Proof.} See Theorem 2 in \cite{Gr}. \hfill $\square$

\vskip 3mm
For any $m,n\in \BZ^+$ with $m<n,$ we define
\begin{equation*}
  \xi_{m,n}:\G^m\lrt \G^n
\end{equation*}
by
\begin{equation*}
  \xi_{m,n}(\gamma):=
  \begin{pmatrix}
    \gamma & 0 \\
    0 & I_{n-m}
  \end{pmatrix},\qquad \gamma\in \G^m.
\end{equation*}
We let
\begin{equation*}
  \G^\infty:=\varinjlim_n \G^n
\end{equation*}
be the inductive limit of the directed system $(\G^n,\xi_{m,n}).$

\begin{definition}\label{def:4.21}
A collection $(f_n)_{n\geq 1}$ is said to be a {\sf stable\
automorphic form} for $\G^\infty$ if it satisfies the following conditions
\eqref{(4.20)} and \eqref{(4.21)}\,:
\begin{equation}\label{(4.20)}
  f_n\in {\bf A}(\G^n),\quad n\geq 1
\end{equation}
and
\begin{equation}\label{(4.21)}
  \mathfrak L_{n+1}f_{n+1}=f_n,\quad n\geq 1.
\end{equation}
\end{definition}

\vskip 3mm
Let
\begin{equation*}
  {\bf A}^\infty={\bf A}(\G^\infty):=\varprojlim_n {\bf A}(\G^n)
\end{equation*}
be the inverse limit of the inverse system $({\bf A}(\G^n),\mathfrak L_{n})$, that is,
the space of all stable automorphic forms for $\G^\infty$.

\vskip 5mm
We propose the following problems.

\begin{problem}\label{prob:4.22}
Discuss the injectivity, the surjectivity and the bijectivity of $\mathfrak L_n.$
\end{problem}

\begin{problem}\label{prob:4.23}
Give examples of stable automorphic forms for
$\G^\infty$.
\end{problem}

\begin{problem}\label{prob:4.24}
Investigate the structure of ${\bf A}^\infty$.
\end{problem}

\begin{remark}\label{rk:4.25}
We refer to \cite{Y1} for more detail on stable automorphic forms for the general linear group
$GL(\infty,\BR)$. We note that the general linear group $GL(n,\BR)\,(n\geq 1)$ is not semisimple.
\end{remark}

\end{section}

\vskip 10mm


\begin{section}
{{\large\bf Applications\ of the stability to\ geometry}}
\setcounter{equation}{0}

\vskip 2mm
In the final section, we give  applications of stable automorphic forms to geometry.

\vskip 3mm
\noindent
{\large\bf 5.1.\ The universal moduli space of abelian varieties}

\vskip 2mm
First of all, for any two nonnegative integers $k,l\in \BZ_+$ with
$k<l,$ we define the mapping $\varphi_{kl}:\BH_k\lrt \BH_l$ by
\begin{equation}\label{(5.1)}
\varphi_{kl}(Z):=\,\begin{pmatrix} Z & 0 \\ 0 & iI_{l-k}
\end{pmatrix},\quad Z\in \BH_k.
\end{equation}
Then the image $\varphi_{kl}(\BH_k)$ is a totally geodesic subspace of $\BH_l$. We let
\begin{equation}\label{(5.2)}
\BH_\infty=\,\varinjlim_k \BH_k
\end{equation}
be the inductive limit of the direct system $(\BH_k,\,\varphi_{kl}).$ $\BH_\infty$ can be described explicitly as follows:
\begin{equation*}
  \left\{ \begin{pmatrix}
            Z & 0 \\
            0 & iI_\infty
          \end{pmatrix} \bigm|\ \, Z\in \BH_k
          \ {\rm for\ some}\ k\geq 1 \right\}.
\end{equation*}
We can show that $\BH_\infty$ is an infinite
dimensional smooth Hermitian symmetric manifold locally closed on $\BC^\infty$, the complex vector space of finite sequences with the finite topology\,(cf.\,\cite{Gl, H}). $\BH_\infty$ has an invariant Riemannian metric which induces the normalized Riemannian metric on each embedded interior subspace $\BH_k$ in $\BH_\infty$.

\vskip 2mm
For each $n\in \BZ^+,$ we put
\begin{equation*}
  G_n:=\,Sp(2n,\BR),\quad K_n:=\,U(n)\quad {\rm and}\quad \G_n:=\,Sp(2n,\BZ).
\end{equation*}
For any $k,l\in \BZ_+$
with $k<l,$ we define the mapping $\pi_{kl}:G_k\lrt G_l$ by
\begin{equation}\label{(5.3)}
\pi_{kl}\left(\begin{pmatrix} A & B\\ C & D\end{pmatrix}\right):=\,
\begin{pmatrix} A & 0 & B & 0 \\ 0 & I_{l-k} & 0 & 0 \\
C & 0 & D & 0 \\ 0 & 0 & 0 & I_{l-k}\end{pmatrix}\qquad {\rm for\ all}\
\begin{pmatrix} A & B\\ C & D\end{pmatrix}\in G_k
\end{equation}
and also define the mapping $\rho_{kl}:\G_k\lrt \G_l$ by the formula \eqref{(5.3)} with
$\begin{pmatrix} A & B\\ C & D\end{pmatrix} \in \G_k.$
Let
\begin{equation*}
  G_\infty:=\varinjlim_k G_k\qquad {\rm and}\qquad
  \G_\infty:=\varinjlim_k \G_k
\end{equation*}
be the inductive limit of the directed systems $(G_k,\pi_{kl})$ and
$(\G_k,\rho_{kl})$ respectively. Then $G_\infty$ and $\G_\infty$
can be described explicitly as follows:
\begin{equation*}
 G_\infty= \left\{ \begin{pmatrix}
             A & 0 & B & 0 \\
             0 & I_\infty & 0 & 0 \\
             C & 0 & D & 0 \\
             0 & 0 & 0 & I_\infty
           \end{pmatrix}\,\Bigm| \ \
           \begin{pmatrix}
             A & B \\
             C & D
           \end{pmatrix}\in G_k\ {\rm for\ some}\ k\geq 1 \right\}
\end{equation*}
and
\begin{equation*}
\G_\infty=\left\{ \begin{pmatrix}
             A & 0 & B & 0 \\
             0 & I_\infty & 0 & 0 \\
             C & 0 & D & 0 \\
             0 & 0 & 0 & I_\infty
           \end{pmatrix}\,\Bigm| \ \
           \begin{pmatrix}
             A & B \\
             C & D
           \end{pmatrix}\in \G_k\ {\rm for\ some}\ k\geq 1 \right\}.
\end{equation*}
We recall that for any two positive integers $k,l\in \BZ^+$ with $k<l,$ the mapping
$u_{k,l}:U(k)\lrt U(l)$
defined by
\begin{equation}\label{(5.4)}
u_{k,l}(A+iB):=\,\begin{pmatrix} A+iB & 0\\ 0 & I_{l-k}\end{pmatrix}\quad
{\rm for\ all}\ A+iB\in U(k)\ {\rm with}\ A,B\in \BR^{(n,n)}
\end{equation}
yields the inductive limit $K_\infty:=U(\infty)$ of the directed system $(U(k),u_{k,l}).$

\begin{lemma}\label{lem:5.1}
Let $k$ and $l$ be two positive integers with
$k<l.$ Then for any $\g\in \G_k$ and $Z\in \BH_k,$ we have
\begin{equation}\label{(5.5)}
\varphi_{kl}(\g\cdot Z)\,=\,\rho_{kl}(\g)\cdot \varphi_{kl}(Z).
\end{equation}
\end{lemma}
\noindent
{\it Proof.} \eqref{(5.5)} follows from an easy computation.
\hfill $\square$

\vskip 3mm
\def\BA{{\mathcal A}}
For each positive integer $k\in \BZ^+,$ we let $\BA_k:=\,\G_k\ba \BH_k$ be
the Siegel modular variety of degree $k$.  We put $\BA_0:=\{ \infty\}.$
According to Lemma \ref{lem:5.1},
for any $k,l\in \BZ^+$ with $k<l,$
we obtain the canonical embedding $s_{kl}:\BA_k\lrt \BA_l$ defined by
\begin{equation}\label{(5.6)}
s_{kl}([Z]):=\,[\varphi_{kl}(Z)]\,=\,\left[
\begin{pmatrix} Z & 0 \\ 0 & iI_{l-k}\end{pmatrix}\right],
\end{equation}
where $[Z]\in \BA_k$ with $Z\in \BH_k$ and $[Z]$ denotes the equivalence class of $Z$.
We let
\begin{equation}\label{(5.7)}
\BA_\infty:=\,\varinjlim_k \BA_k
\end{equation}
be the inductive limit of the directed system $(\BA_k,\,s_{kl}).$
\begin{proposition}\label{prop:5.2}
$G_\infty$ acts on $\BH_\infty$ transitively
and $\G_\infty$ acts on $\BH_\infty$ properly discontinuously.
$\BH_\infty$ is isomorphic to $G_\infty/K_\infty.$ And we have
$$\BA_\infty\,=\,\G_\infty\backslash G_\infty/K_\infty.$$
$\BA_\infty$ is an infinite dimensional Hermitian locally symmetric space.
Furthermore $\BA_\infty$ has a canonical stratification induced from
the canonical stratification of the subspaces
$\mathcal A_{k+1}\setminus\mathcal A_{k}$\ ($\mathcal A_{k+1}$\ setminus\ $\mathcal A_{k}$)
\ with $k\geq 1.$
\end{proposition}
\noindent
{\it Proof.} We observe that $\G_\infty$ is not a finitely generated group. It is countable and an arithmetic discrete subgroup of $G_\infty$. We see that $\G_\infty$ acts on $\BH_\infty$ properly discontinuously and holomorphically. The quotient space is Hausdorff. We can show without difficulty that
\begin{equation*}
 \G_\infty\backslash
 G_\infty/K_\infty =\,\varinjlim_k \BA_k=\mathcal A_\infty.
\end{equation*}
\hfill $\square$

\vskip 2mm
For each nonnegative integer $d\in \BZ_+$, we let $[\G_n,d]$ be the vector space of all
Siegel modular forms on $\BH_n$ of weight $d$. We review some
properties of the Siegel operator $\Phi_{n-1,n}:[\G_n,d]\lrt
[\G_{n-1},d]$\,(\,cf.\,Formula (4.8)\,). According to the theory of singular
modular forms in \cite{Fr1} and \cite{R}, $\Phi_{n-1,n}$ is injective if
$n>2d$ and $\Phi_{n-1,n}$ is an isomorphism if $n>2d+1.$
H. Maass \cite{M} proved that $\Phi_{n-1,n}$ is an isomorphism if
$d$ is even and $d>2n.$

\vskip 3mm
For each nonnegative integer $n\in \BZ_+,$ we put
\begin{equation}\label{(5.8)}
\mathbb A_n:=\,\bigoplus_{d=0}^{\infty}[\G_n,d]\quad \text{for} \ n\geq 1
\quad \text{and} \quad \mathbb A_0:=\BC.
\end{equation}
Then $\mathbb A_n$ is a $\BZ_+$-graded ring which is integrally closed and
of finite type over $\BC:=[\G_n,0].$ We observe that for $m<n,$
the Siegel operator $\Phi_{m,n}$ maps $\mathbb M_n$ into $\mathbb M_m$ preserving
the weights and that $\Phi_{m,n}$ is a ring homomorphism of $\mathbb A_n$
into $\mathbb A_m$. Thus $(\,\mathbb A_n,\,\Phi_{m,n}\,)$ forms an inverse system of
rings over $\BZ_+.$ We let
\begin{equation*}
\mathbb A_\infty:=\,\varprojlim_n \mathbb A_n
\end{equation*}
be the inverse limit of the system $(\mathbb A_n,\,\Phi_{m,n})$. That is,
\begin{equation*}
\mathbb A_\infty\,=\,\left\{\,(f_k)\in \prod_{l\in\BZ_+}
\mathbb A_l\,\Big|\
\Phi_{k,l}(f_l)\,=\,f_k\ \ \ \text{for\ any}\ k<l\ \right\}.
\end{equation*}
If $f=(f_n)\in \mathbb A_\infty,$ then for each $n\in \BZ_+,$ we write
\begin{equation*}
f_n\,=\,\sum_{d=0}^{\infty}\,f_{n,d},\ \ \ \ f_{n,d}\in
[\G_n,d].
\end{equation*}
We note that $\Phi_{m,n}(f_{n,d})\,=\,f_{m,d}$ for all $m,n\in
\BZ_+$ with $m<n.$ For each $d\in \BZ_+,$ the sequence
$\left\{\,(f_{k,d})\,\vert\ k\in \BZ_+\,\right\}$ is called a
{\sf stable\ sequence} of weight $d$. We denote by $S_d$ the
complex vector space consisting of all stable sequences of weight
$d$. Then it is easy to see that
\begin{equation*}
\mathbb A_\infty\,=\,\bigoplus_{d=0}^{\infty} S_d.
\end{equation*}
Then $\mathbb A_\infty$ is a $\BZ_+$-graded ring. It is known that
$\text{dim}_{\BC}\,S_d\,=\,\text{dim}_{\BC}\,[\G_n,d]$ if
$n>2d$\,(\,cf.\,\cite{Fr2},\,p,\,203\,).

\vskip 3mm
Let $S$ be a positive definite even unimodular integral
matrix of degree $m$. Then we define the theta series
$\vartheta_S^{(n)}(Z)$ on $\BH_n$ by
\begin{equation*}
\vartheta_S^{(n)}(Z):=\,\sum_{U\in \BZ^{(m,n)}}\,
e^{\pi i\,{\rm Tr}(S[U]Z)},\ \ \ \ Z\in \BH_n.
\end{equation*}
Here $S[U]:=\,{}^tUSU$ (Siegel's notation).
Then we can show that $\vartheta_S^{(n)}(Z)$ is a Siegel modular form
of weight $m/2$ on $\G_n.$

\vskip3mm
We state the result obtained by Freitag\,\cite{Fr2}.
\begin{theorem}\label{thm:5.3}
$\mathbb A_\infty$ is a polynomial ring in a countably infinite set of indeterminates over $\BC$ given by
\begin{equation*}
\mathbb A_\infty\,=\,\BC[\vartheta_S^{(n)}(Z)\,\vert\ n\in \BZ_+],
\end{equation*}
where $S$ runs over the set of all equivalence classes of irreducible
positive definite symmetric, unimodular even integral matrices.
\end{theorem}
\noindent
{\it Proof.}  See Theorem 2.5 in \cite{Fr2}. \hfill $\square$

\begin{remark}\label{rk:5.4}
The homogeneous quotient field $Q(\mathbb A_\infty)$ of
$\mathbb A_\infty$ is a rational function field with countably infinitely
many variables. But in general $\mathbb A_n$ is not a polynomial ring.
It is well known that the homogeneous function field
$Q(\mathbb A_n)$ of $\mathbb A_n$ is an algebraic function field with the transcendence degree ${\frac 12}n(n+1).$
\end{remark}

\vskip 2mm
For any two nonnegative integers $m,n\in \BZ_+$ with $m<n,$ the Siegel operator $\Phi_{m,n}:
\mathbb A_n\lrt \mathbb A_m$ induces the morphism
$\Phi^{\ast}_{m,n}:\text{Proj}\,\mathbb A_m\lrt
\text{Proj}\,\mathbb A_n$ of projective schemes. The Satake compactification
$\BA_n^{\ast}\,=\,\text{Proj}\,\mathbb A_n$ of $\BA_n$ contains $\BA_n$ as a Zariski open dense subset. As a set, $\BA_n^{\ast}$ is the disjoint union of $\BA_n$ and its rational boundary components, i.e.,
\begin{equation*}
\BA_n^{\ast}\,=\,\BA_n\cup \BA_{n-1}\cup\cdots\BA_1\cup
\BA_0,\ \ \ \BA_0=\left\{{\infty}\right\}.
\end{equation*}
We refer to Satake's paper \cite{Sa}. W. Baily\,\cite{B} proved that $\BA_n^{\ast}$ is a normal projective variety. Obviously $(\,\BA_n^{\ast},\,\Phi_{m,n}^{\ast}\,)$ forms an inductive system of schemes over $\BZ_+.$ We let
\begin{equation*}
\BA^{\ast}_\infty:=\,\lim_{\lrt}\,\BA_n^{\ast}
\end{equation*}
be the inductive limit of $(\,\BA_n^{\ast},\,\Phi_{m,n}^{\ast}\,).$
We call the infinite dimensional variety $\BA^{\ast}_\infty$ the
{\sf universal (or stable)\ Satake\ compactification.}

\vskip 2mm
\begin{theorem}\label{thm:5.5}
The universal Satake compactification $\BA_\infty^{\ast}$ has the following properties:
\vskip 2mm
(1)\ \ $\BA^{\ast}_\infty\,=\,{\rm Proj}\,\,\mathbb A_\infty.$
\vskip 2mm
(2)\ \ $\BA^{\ast}_\infty$ is an infinite dimensional projective
variety which contains $\BA_\infty$ as a Zariski \\
\indent \ \ \ \ \ \ open dense subset.
So $\BA^{\ast}_\infty$ is also called the {\sf Satake\ compactification} of $\BA_\infty.$
\end{theorem}
\noindent
{\it Proof.} The proofs of (1) and (2) follows from Theorem \ref{thm:5.3} and
the following facts (a)--(c)\,:
\vskip 1mm (a) $\mathcal A_n^{\ast}=\,{\rm Proj}\,\,\mathbb A_n$ as schemes
(see \eqref{(5.8)}).
\vskip 1mm (b) For sufficiently large $m,n\in\BZ^+$ with $n>m>2d+1>0$, the Siegel operator
\begin{equation*}
  \Phi_{m,n}:[\G_n,d]\lrt [\G_m,d]
\end{equation*}
\indent \ \ \ \ is an isomorphism.
\vskip 1mm
(c) $\mathcal A_n$ is a Zariski open dense subset of $\BA^{\ast}_n.$
\hfill $\square$

\vskip 5mm
Now we shall describe the analytic local ring of the image of the
boundary point in $\BA_n^{\ast}$ under $f_n^{\ast}$, where
$f_n^{\ast}:\BA_n^{\ast}\lrt \BA^{\ast}_\infty\,(\,n\in \BZ^+\,)$
is the canonical morphism. Let $[Z_k]\in \BA_k\,(\,0\leq k\leq n-1,\
Z_k\in \BH_k\,)$ be a boundary point in
$\BA_n^{\ast}\setminus\BA_n$ ($\BA_n^{\ast}$ setminus $\BA_n$).
We set $Z_{k,\infty}^{\ast}:=\,f_k^{\ast}([Z_k]).$

\vskip 2mm
\begin{theorem}\label{thm:5.6}
The analytic local ring at $Z_{k,\infty}^{\ast}$
in $\BA^{\ast}_\infty$ consists of all sequences
$(f_m)_{m=0}^{\infty}$ with $\Phi_{m,m+1}f_{m+1}=f_m$ such that
each $f_{k+m}\,(m\geq 1)$ is a convergent series of type
\begin{equation*}
f_{k+m}(Z,W_m,\Omega_m)=\sum_{T_m}\phi_{T_m}(Z,W_m)
e^{2\pi i\,{\rm Tr}(T_m\Omega_m)},
\end{equation*}
where $Z$ is an element in a sufficiently small open neighborhood $V$ of
$Z_k$ in $\BH_k$ invariant under the action of
$\G_k,\ W_m\in \BC^{(k,m)},\ \Omega_m\in \BH_m$ and $T_m$ runs over
the set of all semi-positive symmetric half-integral matrices of degree $m$. In addition, each $\phi_{T_m}(Z,W_m)\,(m\geq 1)$ is a Jacobi form
of weight $0$ and index $T_m$ defined on $V\times \BC^{(k,m)}.$
\end{theorem}
\noindent
{\it Proof.}  The proof can be found in \cite{I1}. \hfill $\square$

\vskip 3mm
Ji and Jost \cite{JJ} describe $\mathcal A_\infty^{\ast}$ in a somewhat different way. Since $\BH_k$ is a Hermitian symmetric space of noncompact type, it can be embedded into its compact dual ${\mathfrak Y}_k$ which is a complex projective variety via the Borel embedding. The description of the compact dual ${\mathfrak Y}_k$ is given as follows. We suppose that
$\Lambda=(\BZ^{2k},\langle\ ,\ \rangle)$ is a symplectic lattice with a
symplectic form $\langle\ ,\ \rangle.$ We extend scalars of the lattice $\Lambda$ to $\BC$. Let
\begin{equation*}
{\mathfrak Y}_k:=\left\{\,L\subset \BC^{2k}\,|\ \dim_\BC L=k,\ \
\langle x,y \rangle=0\quad \textrm{for all}\ x,y\in L\,\right\}
\end{equation*}
be the complex Lagrangian Grassmannian variety parameterizing
totally isotropic subspaces of complex dimension $k$. For the
present time being, for brevity, we put $G=Sp(2k,\BR)$ and
$K=U(k).$ The complexification $G_\BC=Sp(2k,\BC)$ of $G$ acts on
${\mathfrak Y}_k$ transitively. If $H$ is the isotropy subgroup of
$G_\BC$ fixing the first summand $\BC^k$, we can identify
${\mathfrak Y}_k$ with the compact homogeneous space $G_\BC/H.$ We
let
\begin{equation*}
{\mathfrak Y}_k^+:=\big\{\,L\in {\mathfrak Y}_k\,|\ -i \langle x,{\bar x}\rangle >0\quad \textrm{for all}\ x(\neq 0)\in L\,\big\}
\end{equation*}
be an open subset of ${\mathfrak Y}_k$. We see that $G$ acts on
${\mathfrak Y}_k^+$ transitively. It can be shown that ${\mathfrak Y}_k^+$ is biholomorphic to $G/K\cong \BH_k.$ A basis of a lattice
$L\in {\mathfrak Y}_k^+$ is given by a unique $2k\times k$ matrix
${}^t(-I_k,\tau)\in \BC^{(2k,k)}$ with $\tau\in\BH_k$. Therefore we can identify
$L$ with $\tau$ in $\BH_k$. In this way, we embed $\BH_k$ into
${\mathfrak Y}_k$ as an open subset of ${\mathfrak Y}_k$.

\vskip 2mm
The closure $\overline{\BH_k}$ of $\BH_k$ in ${\mathfrak Y}_k$ is compact. The standard embedding $\varphi_{kl}$ (see Formula \eqref{(5.1)}) of $\BH_k$ into the boundary of $\overline{\BH_l}$
with $k< l$ and the translates by $\G_l:=Sp(2l,\BZ)$ of these standard boundary components
give all the rational boundary components (briefly
${\bf rbc}$) of $\BH_k$. We denote by ${\overline{\BH_k}}_{,\BQ}$
the union of $\BH_k$ with these ${\bf rbc}$. Then there exists the so-called Satake topology on ${\overline{\BH_k}}_{,\BQ}$ such that $\G_k$ acts continuously on ${\overline{\BH_k}}_{,\BQ}$. Then we obtain the Satake compactification of $\BH_k$

\begin{equation*}
  \mathcal A_k^{\ast}=\G_k\backslash {\overline{\BH_k}}_{,\BQ}.
\end{equation*}

\vskip 2mm
From the increasing sequence
\begin{equation*}
  {\overline{\BH_1}}_{,\BQ} \hookrightarrow {\overline{\BH_2}}_{,\BQ} \hookrightarrow   {\overline{\BH_3}}_{,\BQ} \hookrightarrow
  \cdots \hookrightarrow  {\overline{\BH_k}}_{,\BQ} \cdots,
\end{equation*}
we get the inductive limit
\begin{equation*}
  {\overline{\BH_\infty}}_{,\BQ}=\varinjlim_k {\overline{\BH_k}}_{,\BQ}.
\end{equation*}

${\overline{\BH_\infty}}_{,\BQ}$ can be realized as
\begin{equation*}
  {\overline{\BH_\infty}}_{,\BQ}=
  \left\{ \begin{pmatrix}
            \Omega & 0 \\
            0 & iI_\infty
            \end{pmatrix}
            \Bigg| \ \ \Omega \in {\overline{\BH_k}}_{,\BQ}\ {\rm for\ some}\ k
            \right\}.
\end{equation*}

Taking the quotient of ${\overline{\BH_\infty}}_{,\BQ}$ by $\G_\infty$, we obtain the completion of $\mathcal A_\infty$,
\begin{equation*}
 \mathcal A_\infty^{\rm Sat}= \G_\infty \backslash {\overline{\BH_\infty}}_{,\BQ}.
\end{equation*}

Since
\begin{equation*}
 \mathcal A_\infty^{\rm Sat}= \varinjlim_k \mathcal A_k^{\ast}=
 \bigcup_{k\geq 0} \mathcal A_k^{\ast}
\end{equation*}
under the inclusion $\mathcal A_k^{\ast}\hookrightarrow
\mathcal A_{k+1}^{\ast}$, we see that
$\mathcal A_\infty^{\rm Sat}=\mathcal A_\infty^{\ast}$
(cf. Theorem 5.2).

\vskip 3mm
Ji and Jost \cite{JJ} obtain the following result.
\begin{proposition}\label{prop:5.7}
The universal Satake compactification $\mathcal A_\infty^{\ast}$
admits the following decomposition
 \begin{equation*}
 \mathcal A_\infty^{\ast}=\mathcal A_\infty \bigsqcup
 (\mathcal A_0 \sqcup \mathcal A_1 \sqcup \mathcal A_2 \sqcup
 \cdots ),
 \end{equation*}
 where
 \begin{equation*}
 \bigsqcup_{k\geq 0} \mathcal A_k
 \end{equation*}
 is the boundary, and $\mathcal A_\infty$ is the interior in some sense, which can also be decomposed into a non-disjoint union of
 $\mathcal A_k,\ k\geq 0$. Every $\mathcal A_k$ can appear in
 $\mathcal A_\infty^{\ast}$ in two ways: either in the interior
 $\mathcal A_\infty$ or in the boundary
 $\bigsqcup_{k\geq 0} \mathcal A_k$.
\end{proposition}
\noindent
{\it Proof.}  The proof can be found in section 3 of the paper \cite{JJ} of Ji and Jost. \hfill $\square$

\vskip 7mm
\noindent
{\large\bf 5.2.\ The universal moduli space of curves}

\vskip 2mm
For each positive integer $g\in \BZ^+$, we let $\mathcal M_g$ be the moduli space of projective curves of genus $g$ and $\mathcal A_g=Sp(2g,\BZ)\backslash \BH_g$ the Siegel modular variety of degree $g$.
According to Torelli's theorem, the Jacobi mapping
\begin{equation}\label{(5.9)}
T_g:{\Cal M}_g \lrt {\Cal A}_g
\end{equation}
defined by
\begin{equation*}
C \longmapsto J(C):= {\rm the\ Jacobian\ of}\ C
\end{equation*}
is injective, and in fact it is an embedding. $T_g$ induces an embedding
\begin{equation*}
T_g^{\ast}:{\mathcal M}_g^{\ast} \lrt {\mathcal A}_g^{\ast},
\end{equation*}
where ${\mathcal M}_g^{\ast}$ (resp.\,${\mathcal A}_g^{\ast}$) is the Satake compactification of ${\mathcal M}_g$
(resp.\,${\mathcal A}_g$). The Jacobian locus
$J_g:=T_g({\mathcal M}_g)$ is a $(3g-3)$-dimensional subvariety of ${\mathcal A}_g$ if $g\geq 2.$
Let $J_g^{\ast}:=T_g^{\ast}({\mathcal M}_g^{\ast})$ be the Satake
compactification of $J_g$, which is equal to the closure of $J_g$ in
${\mathcal A}_g^{\ast}$ for the Satake topology.

\vskip 2mm
For convenience, we set
\begin{equation*}
  J_0= J_0^{\ast}=  {\mathcal A}_0={\mathcal A}_0^{\ast}=\{\infty\},\quad {\rm one\ point}.
\end{equation*}
We define
\begin{equation*}
  J_\infty :=\bigcup_{g\geq 0} J_g
\end{equation*}
and
\begin{equation*}
{\mathcal A}_\infty=
  \bigcup_{g\geq 0} {\mathcal A}_g,\qquad
{\mathcal A}_\infty^{\ast}=
  \bigcup_{g\geq 0} {\mathcal A}_g^{\ast}=\varinjlim_g
   {\mathcal A}_g^{\ast} \qquad ({\rm see}\ (5.7)\ {\rm and} \ (5.13)).
\end{equation*}

\begin{proposition}\label{prop:5.8}
(1) The boundary of $J_g^{\ast}$ is the union of $J_{g_1}\times \cdots \times J_{g_k}$,\ where $g_1+\cdots+g_k\leq g$ with
$k\geq 1.$
\vskip 2mm\noindent
(2) For any two positive integers $k,g$ with $k<g$, if $J_k$ appears in the boundary of $J_g^{\ast}$, then the closure of $J_k$ is equal to the Satake compactification $J_k^{\ast}$ of $J_k$.
\vskip 2mm\noindent
(3) The subspace $J_\infty^{\ast}$ of ${\mathcal A}_\infty^{\ast}$
has a canonical stratification such that the closure of each stratum is a projective variety over $\BC$, and $J_\infty^{\ast}$ is the Satake compactification of $J_\infty$.
$J_g^{\ast}$ can appear in many different ways in $J_\infty^{\ast}$.
\vskip 2mm\noindent
(4) $J_\infty^{\ast}$ is connected in ${\mathcal A}_\infty^{\ast}$.
\vskip 2mm\noindent
(5) For any $g\in\BZ^+,$ there is a unique way to embed $J_g$ into $J_{g+1}^{\ast}$ which is the closure of $J_{g+1}$ inside
$J_\infty^{\ast}$. Under this inclusion, we get an increasing sequence of spaces
\begin{equation*}
  J_0^\ast \hookrightarrow J_1^\ast \hookrightarrow J_2^\ast
  \hookrightarrow J_3^\ast \hookrightarrow \cdots
\end{equation*}
and
\begin{equation*}
J_\infty^{\ast} :=\bigcup_{g\geq 0} J_g^{\ast}=\varinjlim_g J_g^{\ast}.
\end{equation*}
\end{proposition}
\vskip 2mm\noindent
\begin{proof} The proof can be found in section 4 of the paper \cite{JJ} of Ji and Jost.
\end{proof}

\begin{theorem}\label{thm:5.9}
For any $g\in\BZ^+$, there exists a Riemannian metric on $J_g^{\ast}$ that induces
a Riemannian metric on each stratum.
As a result, there exists a measure on $J_g^{\ast}$ that induces
a finite volume measure on each stratum.
\end{theorem}
\begin{proof}
The proof can be found in \cite[Theorem 5.2 and Corollary 5.3]{JJ}.
\end{proof}

\vskip 2mm
\begin{definition}\label{def:5.10}
A modular form $f\in [\G_g,k]$ is called a {\sf Schottky-Siegel form} of weight $k$ for $J_g$ (resp. $Hyp_g$) if it vanishes along
$J_g$ (resp. $Hyp_g$). A collection $(f_g)_{g\geq 0}$ is called a {\sf stable Schottky-Siegel form} of weight $k$ for the Jacobian locus (resp. the hyperelliptic locus) if $(f_g)_{g\geq 0}$ is a stable modular form of weight $k$ and $f_g$ vanishes along $J_g$ (resp. $Hyp_g$) for every $g\geq 0.$
\end{definition}

G.~Codogni and N.~I. Shepherd-Barron \cite{Cod-Sh} proved the following.
\begin{theorem}\label{thm:5.11}
There do not exist stable Schottky-Siegel form for the Jacobian locus.
\end{theorem}
\begin{proof}
See \cite[Theorem 1.3 and Corollary 1.4]{Cod-Sh}.
\end{proof}

\begin{remark}\label{rk:5.12}
Let
\begin{equation}\label{(5.10)}
  \varphi_g(\tau):=\theta_{E_8\oplus E_8,g}(\tau)
  -\theta_{D_{16}^+,g}(\tau),\qquad \tau\in\BH_g
\end{equation}
be the Igusa modular form, that is, the difference of the theta series in genus $g$ associated to the two distinct positive even unimodular quadratic forms $E_8\oplus E_8$ and $D_{16}^+$ of rank $16$. We see that $\varphi_g(\tau)$ is a Siegel modular form on $\BH_g$ of weight $8$. Since $\Phi_{g-1,g} \varphi_g=\varphi_{g-1}$ for all $g\geq 1$, a collection $(\varphi_g)_{g\geq 0}$ is a stable
modular form of weight $8$. Igusa \cite{I2,I3} showed that the Schottky-Siegel form discovered by Schottky \cite{Sch} is an explicit rational multiple of $\varphi_4$. In \cite{I2}, he also
showed that the Jacobian locus $J_4$ is reduced and irreducible, and so cuts out exactly $J_4$ in $\mathcal A_4.$ Indeed, $\varphi_4(\tau)$ is a degree $16$ polynomial in the Thetanullwerte of genus $4$. On the other hand, Grushevsky and Salvati Manni \cite{G-SM} showed that the Igusa modular form $\varphi_5$ of genus $5$ cuts out exactly the trigonal locus in $J_5$ and so does not vanish along $J_5$. Thus $(\varphi_g)_{g\geq 0}$ is not a stable
Schottky-Siegel form.
\end{remark}

\vskip 2mm
G.~Codogni \cite{Cod} proved the following.
\begin{theorem}\label{thm:5.13}
There exist non-trivial stable Schottky-Siegel form for the hyperelliptic locus. Precisely the ideal of stable Schottky-Siegel forms for the hyperelliptic locus is generated by differences of theta series
\begin{equation*}
\Theta_P-\Theta_Q,
\end{equation*}
where $P$ and $Q$ are positive definite even unimodular quadratic forms of the same rank.
See Definition \ref{def:4.10} for the definition of $\Theta_P$.
\end{theorem}
\begin{proof} The proof can be found in \cite[Theorem 1.2]{Cod}.
\end{proof}

\vskip 1mm

\begin{remark}\label{rk:5.14}
Let $\varphi_g(\tau)$ be the Igusa modular form defined by the formula \eqref{(5.10)}.
We denote by $[\G_g,k]_0$ be the space of all Siegel cuspidal Hecke eigenforms on $\BH_g$ of weight $k$.
It is known that $[\G_4,8]_0=\BC\cdot\varphi_4$ (for a nice proof of this, we refer to \cite{DI}).
Poor \cite{P} showed that $\varphi_g(\tau)$ vanishes on the hyperelliptic locus $Hyp_g$ for all $g\geq 1$, and the divisor of $\varphi_g(\tau)$ in $\mathcal{A}_g$ is proper and irreducible for all $g\geq 4$.
And Ikeda \cite{Ik}
proved that if $g\equiv k\,({\rm mod}\,2),$ there exists a canonical lifting
\begin{equation*}
  I_{g,k}:[\G_1,2k]_0\lrt [\G_{2g},g+k]_0.
\end{equation*}
Considering the special cases of the Ikeda lift maps $I_{2,6}$ and $I_{6,6}$, Breulman and Kuss \cite{BK} showed that
\begin{equation*}
  I_{2,6}(\Delta)=c\,\varphi_4, \quad c(\neq 0)\in\BC,
\end{equation*}
and constructed a nonzero Siegel cusp form of degree $12$ and weight $12$ which is the image of $\Delta(\tau)$ under the lifting $I_{6,6}$, where
\begin{equation*}
\Delta(\tau)=(2\pi i)^{12}~q\prod_{n=1}^{\infty}
\left( 1-q^n\right)^{24},\quad q:=e^{2\pi i\tau},\quad \tau\in\BH_1
\end{equation*}
is a cusp form of weight 12.
\end{remark}

\vskip 7mm
\noindent
{\large\bf 5.3.\ The universal moduli space of polarized real tori}
\newcommand{\La}{\Lambda}
\newcommand{\FA}{\mathfrak A}
\newcommand{\FL}{\mathfrak L}

\vskip 3mm
Let
\begin{equation*}
  \CP:=\left\{ Y\in\BR^{(n,n)}\,\vert\ Y=\,{}^tY >0\,\right\}
\end{equation*}
be the cone of positive definite symmetric real matrices of degree $n$. Then $GL(n,\BR)$ acts on $\CP$ transitively by
\begin{equation}\label{(5.11)}
g\circ Y:=gY\,{}^tg,\qquad g\in GL(n,\BR),\ Y\in \CP.
\end{equation}

\vskip 3mm
First we recall the concept of polarized real tori\,(cf.\,
\cite[p.\,295]{Y}).
\begin{definition}\label{def:5.15}
A real torus $T=\BR^n/\Lambda$ with a lattice $\La$ in $\BR^n$ is said to be $\textsf{polarized}$
if the the associated complex torus $\FA=\BC^n/L$ is a polarized real abelian variety, where
$L=\,\BZ^n+\,i\,\La$ is a lattice in $\BC^n.$ Moreover if $\FA$ is a principally polarized real abelian
variety, $T$ is said to be \textsf{principally polarized}.
Let $\Phi:T\lrt\FA$ be the smooth embedding of $T$ into $\FA$ defined by
\begin{equation}\label{(5.12)}
\Phi(v+\La):=\,i\,v\,+\,L, \qquad v\in\BR^n.
\end{equation}
Let $\FL$ be a polarization of $\FA$, that is, an ample line bundle over $\FA$. The pullback $\Phi^*\FL$ is
called a \textsf{polarization} of $T$. We say that a pair $(T,\Phi^*\FL)$ is a \textsf{polarized real torus}.
\end{definition}

\begin{example}\label{exa:5.16}
Let $Y\in \CP$ be a $n\times n$ positive definite symmetric real matrix. Then
$\La_Y=\,Y\BZ^n$ is a lattice in $\BR^n$. Then the $n$-dimensional torus $T_Y=\BR^n/\La_Y$ is a
principally polarized real torus. Indeed,
\begin{equation*}
\FA_Y\,=\,\BC^n/L_Y, \qquad  L_Y\,=\BZ^n+\,i\,\La_Y
\end{equation*}
is a princially polarized real abelian variety.
Its corresponding Hermitian form $H_Y$ is given by
\begin{equation*}
H_Y(x,y)\,=\,E_Y(i\,x,y)\,+\,i\,E_Y(x,y)\,=\,{}^tx\,Y^{-1}\,{\overline y},\qquad x,y\in\BC^n,
\end{equation*}
where
$E_Y$ denotes the imaginary part of $H_Y.$ It is easily checked that $H_Y$ is positive definite and
$E_Y(L_Y\times L_Y)\subset \BZ$\,(cf.\,\cite[pp.\,29--30]{Mum}).
The real structure $\s_Y$ on $\FA_Y$ is a complex conjugation.
In addition, if $\det Y=1$, the real torus $T_Y$ is said to be ${\sf special}$.
We refer to \cite[pp.\,275--279]{Y} for more details about {\sf real structure}.
\end{example}

\begin{example}\label{exa:5.17}
Let $Q=\,\begin{pmatrix} \sqrt{2} & \ \sqrt{3} \\ \sqrt{3} & -\sqrt{5}
\end{pmatrix}$ be a $2\times 2$ symmetric real matrix of signature $(1,1)$.
Then $\La_Q=\,Q\BZ^2$ is a lattice in $\BR^2$. Then the real torus $T_Q=\,\BR^2/\La_Q$ is not
polarized because the associated complex torus $\FA_Q=\,\BC^2/L_Q$ is not an abelian variety, where
$L_Q=\BZ^2+\,i\,\La_Q$ is a lattice in $\BC^2$.
\end{example}

\begin{definition}\label{def:5.18}
Two polarized tori $T_1\,=\,\BR^n/\La_1$ and $T_2\,=\,\BR^n/\La_2$ are said to be isomorphic if the associated
polarized real abelian varieties $\FA_1\,=\,\BC^n/L_1$ and $\FA_2\,=\,\BC^n/L_2$ are isomorphic, where
$L_i\,=\,\BZ^n\,+\,i\,\La_i\ (i=1,2),$ more precisely, if there exists a linear isomorphism $\varphi:\BC^n\lrt \BC^n$ such that
\begin{eqnarray}\label{(5.13)-(5.15)}
\varphi (L_1)&=&L_2,\\
\varphi_* (E_1)&=&E_2,\\
\varphi_* (\s_1) &=& \varphi\circ \s_1\circ \varphi^{-1}\,=\,\s_2,
\end{eqnarray}
where $E_1$ and $E_2$ are polarizations of $\FA_1$ and $\FA_2$ respectively, and $\s_1$ and $\s_2$ denotes
the real structures (in fact complex conjugations) on $\FA_1$ and $\FA_2$ respectively.
\end{definition}

\begin{example}\label{exa:5.19}
Let $Y_1$ and $Y_2$ be two $n\times n$ positive definite symmetric real matrices. Then
$\La_i:=\,Y_i\,\BZ^n$ is a lattice in $\BR^n$ $(i=1,2)$. We let
$$T_i:=\,\BR^n/\La_i,\qquad i=1,2$$
be real tori of dimension $n$. Then according to Example 5.16, $T_1$ and $T_2$ are
principally polarized real tori.
We see that $T_1$ is isomorphic to $T_2$ as polarized real tori if and only if there is
an element $A\in GL(n,\BZ)$ such that $Y_2\,=\,A\,Y_1\,{}^tA.$
\end{example}

\begin{example}\label{exa:5.20}
Let $Y=\,\begin{pmatrix} \sqrt{2} & \ \sqrt{3} \\ \sqrt{3} & \sqrt{5}
\end{pmatrix}$. Let $T_Y\,=\,\BR^2/\La_Y$ be a two dimensional principally polarized torus, where $\La_Y=\,Y\BZ^2$ is a lattice
in $\BR^2.$ Let $T_Q$ be the torus in Example 5.17. Then $T_Y$ is diffeomorphic to $T_Q$.
But $T_Q$ is not polarized.
$T_Y$ admits a differentiable embedding into a complex projective space but $T_Q$ does not.
\end{example}
\vskip 5mm
Let
\begin{equation*}
  G^n=SL(n,\BR),\quad K^n=SO(n)\quad {\rm and}\quad \G^n=GL(n,\BZ)/\{ \pm I_n\}.
\end{equation*}

\vskip 3mm
Let
\begin{equation*}
  {\mathbb X}_n:=\left\{ Y\in\CP \,\vert\ \det Y=1\,\right\}
\end{equation*}
be the subspace of $\CP$. We see that $G^n$ acts on ${\mathbb X}_n$
transitively via \eqref{(5.11)}, and $K^n$ is the stabilizer at $I_n$.
Thus ${\mathbb X}_n$ is a symmetric space which is diffeomorphic to the
homogeneous space $G^n/K^n$ through the following correspondence
\begin{equation*}
G^n/K^n \lrt {\mathbb X}_n,\qquad gK^n\mapsto g\circ I_n=g\,{}^tg,\quad g\in G_n.
\end{equation*}

\vskip 3mm
The arithmetic variety
\begin{equation*}
  \mathfrak R_n=GL(n,\BZ)\backslash \CP
  =GL(n,\BZ)\backslash GL(n,\BR)/O(n)
\end{equation*}
is the moduli space of principally polarized real tori of dimension $n$. Let
\begin{equation}\label{(5.16)}
  \mathfrak X_n:=\G^n\backslash {\mathbb X}_n=\G^n\backslash G^n/ K^n
\end{equation}
be the moduli space of ${\sf special}$ principally polarized real tori of dimension $n$.

\vskip 5mm
For any two positive integers $m,n\in \BZ^+$ with $m<n$, we define
\begin{equation*}
  \eta_{m,n}:G^m\lrt G^n
\end{equation*}
by
\begin{equation}\label{(5.17)}
  \eta_{m,n}(A):=
  \begin{pmatrix}
    A & 0 \\
    0 & I_{n-m}
  \end{pmatrix}\qquad {\rm for\ all}\ A\in G^m.
\end{equation}
We let
\begin{equation*}
  G^\infty:=\varinjlim_n G^n,\qquad K^\infty:=\varinjlim_n K^n
  \quad {\rm and}\quad \G^\infty:=\varinjlim_n \G^n
\end{equation*}
be the inductive limits of the directed systems $(G^n,\eta_{m,n}),
\ (K^n,\eta_{m,n})$ and $(\G^n,\eta_{m,n})$ respectively.

\vskip 3mm
Let $J_n$ (resp.\,${\rm Hyp}_n$) be the Jacobian locus (resp.\, the hyperelliptic locus) in the Siegel modular variety $\mathcal A_n.$
We define
\begin{equation*}
 {\mathbb X}_{n,J}:=\left\{  Y\in {\mathbb X}_n\,\vert\ \mathfrak A_Y\
 {\rm is\ the\ Jacobian\ of\ a\ curve}\ {\rm of\ genus}\ n,
 \ i.e., \ [\mathfrak A_Y]\in J_n \,   \right\}
\end{equation*}
and
\begin{equation*}
 {\mathbb X}_{n,H}:=\left\{  Y\in {\mathbb X}_n\,\vert\ \mathfrak A_Y\
 {\rm is\ the\ Jacobian\ of\ a\ hyperelliptic\ curve}\
 {\rm of\ genus}\ n\,  \right\}.
\end{equation*}
See Example \ref{exa:5.16} for the definition of $\mathfrak A_Y.$
We see that $\G^n$ acts on both ${\mathbb X}_{n,J}$ and ${\mathbb X}_{n,H}$
properly discontinuously. So we may define
\begin{equation*}
  \mathfrak X_{n,J}:=\G^n\backslash {\mathbb X}_{n,J}\qquad {\rm and} \qquad
  \mathfrak X_{n,H}:=\G^n\backslash {\mathbb X}_{n,H}.
\end{equation*}
$\mathfrak X_{n,J}$ and $\mathfrak X_{n,H}$ are called the
${\sf Jacobian\ real\ locus}$ and the ${\sf hyperelliptic\ real\ locus}$ respectively.

\begin{problem}\label{prob:5.21}
{\bf Problem.} Characterize the Jacobian real locus.
\end{problem}
\vskip 3mm
Let $\mathfrak X_n^{S}$ be the standard or maximal Satake compactification of
$\mathfrak X_n.$ For the details of the standard or maximal Satake compactification
of a locally symmetric space, we refer to \cite[pp.\,286--291]{BoJ1}, \cite{BoJ2} and
\cite[pp.\,7--9]{W-W}.
We denote by $\mathfrak X_{n,J}^{S}$ (resp.\,
$\mathfrak X_{n,H}^{S}$) the standard or maximal Satake compactification of
$\mathfrak X_{n,J}$ (resp.\,$\mathfrak X_{n,H}$). We can show that
$\mathfrak X_{n,J}^{S}$ (resp.\,$\mathfrak X_{n,H}^{S}$) is the closure of $\mathfrak X_{n,J}$ (resp.\,$\mathfrak X_{n,H}$) inside
$\mathfrak X_n^{S}$. We have the increasing sequences

\begin{equation*}
  \mathfrak X_1^{S}\hookrightarrow \mathfrak X_2^{S}\hookrightarrow
  \mathfrak X_3^{S}\hookrightarrow \cdots,
\end{equation*}

\begin{equation*}
  \mathfrak X_{1,J}^{S}\hookrightarrow \mathfrak X_{2,J}^{S}\hookrightarrow
  \mathfrak X_{3,J}^{S}\hookrightarrow \cdots
\end{equation*}
and
\begin{equation*}
  \mathfrak X_{1,H}^{S}\hookrightarrow \mathfrak X_{2,H}^{S}\hookrightarrow
  \mathfrak X_{3,H}^{S}\hookrightarrow \cdots.
\end{equation*}

\vskip 5mm\noindent
We put
\begin{equation*}
 \mathfrak X_{\infty}^{S}:=\varinjlim_n \mathfrak X_n^S,\qquad
 \mathfrak X_{\infty,J}^{S}:=\varinjlim_n \mathfrak X_{n,J}^S\quad  {\rm and}\quad
 \mathfrak X_{\infty,H}^{S}:=\varinjlim_n \mathfrak X_{n,H}^S.
\end{equation*}

\vskip 3mm \noindent
For any two positive integers $m,n\in\BZ^+$ with $m<n$, we embed
${\mathbb X}_m$ into ${\mathbb X}_n$ as follows:
\begin{equation*}
 \psi_{m,n}:{\mathbb X}_m \lrt {\mathbb X}_n, \qquad
 Y \mapsto \begin{pmatrix}
              Y & 0 \\
              0 & I_{n-m}
            \end{pmatrix}\qquad {\rm for\ all}\ Y\in {\mathbb X}_m.
\end{equation*}
We let
\begin{equation*}
  {\mathbb X}_\infty= \varinjlim_n {\mathbb X}_n
\end{equation*}
be the inductive limit of the directed system $({\mathbb X}_n,\psi_{m,n}).$
We can show that
\begin{equation*}
  {\mathbb X}_\infty= G^\infty/ K^\infty.
\end{equation*}

\vskip 3mm
Now we have the Grenier operator
\begin{equation*}
  \mathfrak L_n: {\bf A}(\G^n)\lrt {\bf A}(\G^{n-1})
\end{equation*}
defined by the formula \eqref{(4.17)}.

\begin{definition}\label{def:5.22}
An automorphic form $f\in {\bf A}(\G^n)$ is said to be a ${\sf Grenier\!-\!Schottky\ automorphic}$
${\sf form}$ for the Jacobian real locus (resp.\,the hyperelliptic real locus) if it vanishes along
$\mathfrak X_{n,J}$ (resp.\,$\mathfrak X_{n,H}$).
A collection $(f_n)_{n\geq 1}$ is called a ${\sf stable\ Grenier\!-\!Schottky\
automorphic\ form}$ for the Jacobian real locus (resp.\,the hyperelliptic real locus) if it satisfies the following conditions
{\rm (SGS1)} and {\rm (SGS2)}\,:
\vskip 2mm
{\rm (SGS1)} \ \ $f_n$ is a Grenier-Schottky automorphic form
for the Jacobian real locus \\
\indent\indent \ \ \ \ \ \ \ \ \,(resp.\,the hyperelliptic real locus)
for each $n\geq 1.$
\vskip 2mm
{\rm (SGS2)}\ \  $\mathfrak L_n f_n=f_{n-1}$ for all $n> 1.$
\end{definition}

\vskip 3mm
The following natural question arises\,:
\begin{question}\label{que:5.23}
Are there stable Grenier-Schottky automorphic forms for the Jacobian real locus (resp.\,the hyperelliptic real locus) ?
\end{question}

\end{section}

\vskip 5mm


\end{document}